\documentclass[reqno, oneside]{amsart}

\usepackage{amsthm}

\usepackage{
amssymb,amsmath,enumerate,stackrel,tabu,mathrsfs,verbatim,stmaryrd,
leftidx}
\usepackage{xcolor} 
\usepackage{xr-hyper}
\usepackage{nicefrac,stackrel}
\usepackage[all,cmtip]{xy}
\usepackage[utf8]{inputenc} 
\chardef\cprime"7E 
\usepackage{cases}

\usepackage[T1]{fontenc} 

\usepackage[pagebackref,hyperindex,linktocpage=true]{hyperref}
\hypersetup{
    colorlinks,
    linkcolor={red!50!black},
    citecolor={blue!50!black},
    urlcolor={blue!50!black},
    filecolor={red!50!black} 
}

\usepackage{tikz}
\usetikzlibrary{matrix,arrows}

\newcounter{heyheyCounter}[section]

\newif \ifDraft 
\Draftfalse

\ifDraft
  \usepackage[top=1.5cm, bottom=1.5cm, left=.5cm, right=7.5cm, heightrounded, marginparwidth=7cm, marginparsep=.2cm]{geometry}
  \newcommand{\heyhey}[1]{\refstepcounter{heyheyCounter}$\bigstar^\text{\theheyheyCounter}$\marginpar{\footnotesize $\bigstar^\text{\theheyheyCounter}$ #1}}
\else
  \usepackage[top=1.5cm, bottom=1.5cm, left=3cm, right=3cm, heightrounded]{geometry}
  \usepackage{letltxmacro}
  \LetLtxMacro\Oldfootnote\footnote
  \renewcommand{\footnote}[2][]{\relax}
  \newcommand{\heyhey}[1]{}
\fi

\newcommand{\timo}[1]{\heyhey{Timo: #1}}

\usepackage{xcolor}

\usepackage{theoremref}

\makeatletter
\@addtoreset{equation}{section}
\makeatother

\numberwithin{equation}{section}

\theoremstyle{definition}
\newtheorem{Defi}{Definition}[section] \newcommand{\defi}{\begin{Defi}} \newcommand{\xdefi}{\end{Defi}} 
\newtheorem{Cons}[Defi]{Construction} \newcommand{\cons}{\begin{Cons}} \newcommand{\xcons}{\end{Cons}} 
\newtheorem{DefiLemm}[Defi]{Definition and Lemma} \newcommand{\defilemm}{\begin{DefiLemm}} \newcommand{\xdefilemm}{\end{DefiLemm}} 
\newtheorem{Bsp}[Defi]{Example} \newcommand{\exam}{\begin{Bsp}} \newcommand{\xexam}{\end{Bsp}} 
\newtheorem{Syno}[Defi]{Synopsis} \newcommand{\syno}{\begin{Syno}} \newcommand{\xsyno}{\end{Syno}} 
\newtheorem{Bem}[Defi]{Remark} \newcommand{\rema}{\begin{Bem}} \newcommand{\xrema}{\end{Bem}} 
\newtheorem{Notation}[Defi]{Notation} \newcommand{\nota}{\begin{Notation}} \newcommand{\xnota}{\end{Notation}} 

\theoremstyle{plain}
\newtheorem{Theo}[Defi]{Theorem} \newcommand{\theo}{\begin{Theo}} \newcommand{\xtheo}{\end{Theo}} 
\newtheorem{Satz}[Defi]{Proposition} \newcommand{\prop}{\begin{Satz}} \newcommand{\xprop}{\end{Satz}} 
\newtheorem{Lemm}[Defi]{Lemma} \newcommand{\lemm}{\begin{Lemm}} \newcommand{\xlemm}{\end{Lemm}} 
\newtheorem{Coro}[Defi]{Corollary} \newcommand{\coro}{\begin{Coro}} \newcommand{\xcoro}{\end{Coro}}
\newtheorem{Ques}[Defi]{Question} \newcommand{\ques}{\begin{Ques}} \newcommand{\xques}{\end{Ques}}
\newtheorem{Conj}[Defi]{Conjecture} \newcommand{\conj}{\begin{Conj}} \newcommand{\xconj}{\end{Conj}}

\newcommand{\refit}[1]{(\ref{item--#1})}
\newcommand{\refeq}[1]{(\ref{eqn--#1})}

\newcommand{\eqn}{\begin{equation}} \newcommand{\xeqn}{\end{equation}}
\newcommand{\eqnarr}{\begin{eqnarray*}} \newcommand{\xeqnarr}{\end{eqnarray*}}
\newcommand{\eqnarra}{\begin{eqnarray}} \newcommand{\xeqnarra}{\end{eqnarray}}

\newcommand{\pf}{\begin{proof}} \newcommand{\xpf}{\end{proof}}

\newcommand{\nc}{\newcommand}
\nc{\StP}[1]{\cite[\href{http://stacks.math.columbia.edu/tag/#1}{Tag #1}]{StacksProject}}
\nc{\StPz}[2]{\cite[Tags~\href{http://stacks.math.columbia.edu/tag/#1}{#1},
\href{http://stacks.math.columbia.edu/tag/#2}{#2}]{StacksProject}} 
\nc{\StPd}[3]{\cite[Tags~\href{http://stacks.math.columbia.edu/tag/#1}{#1},
\href{http://stacks.math.columbia.edu/tag/#2}{#2},
\href{http://stacks.math.columbia.edu/tag/#3}{#3}]{StacksProject}} 

\nc{\on}{\operatorname}
\nc{\aff}{{\on{aff}}}
\nc{\modi}{{\on{mod}}} 
\nc{\even}{{\on{even}}}
\nc{\odd}{{\on{odd}}}
\nc{\naive}{{\on{naive}}}
\nc{\hofib}{\on{hofib}}
\nc{\Bun}{\on{Bun}}
\nc{\ad}{{\on{ad}}}
\nc{\lft}{{\on{lft}}}

\nc{\Weil}{{\on{Weil}}} 
\nc{\FWeil}{{\on{FWeil}}} 
\nc{\constr}{{\on{cons}}} 

\nc{\str}{\on{-}}
\nc{\perf}{{\on{perf}}}
\nc{\Rel}{{\on{Pos}}}
\nc{\lan}{\langle}
\nc{\ran}{\rangle}

\nc{\bbA}{{\mathbb A}} 
\nc{\bbB}{{\mathbb B}}
\nc{\bbC}{{\mathbb C}}
\nc{\bbD}{{\mathbb D}}
\nc{\bbE}{{\mathbb E}}
\nc{\bbF}{{\mathbb F}}
\nc{\bbG}{{\mathbf G}}
\nc{\bbH}{{\mathbb H}}
\nc{\bbI}{{\mathbb I}}
\nc{\bbJ}{{\mathbb J}}
\nc{\bbK}{{\mathbb K}}
\nc{\bbL}{{\mathbb L}}
\nc{\bbM}{{\mathbb M}}
\nc{\bbN}{{\N}} 
\nc{\bbO}{{\mathbb O}}
\nc{\bbP}{{\mathbb P}} 
\nc{\bbQ}{{\mathbb Q}} 
\nc{\bbR}{{\mathbb R}}
\nc{\bbS}{{\mathbb S}}
\nc{\bbT}{{\mathbb T}}
\nc{\bbU}{{\mathbb U}}
\nc{\bbV}{{\mathbb V}}
\nc{\bbW}{{\mathbb W}}
\nc{\bbX}{{\mathbb X}}
\nc{\bbY}{{\mathbb Y}}
\nc{\bbZ}{{\mathbb Z}}

\nc{\calA}{{\mathcal A}}
\nc{\calB}{{\mathcal B}}
\nc{\calC}{{\mathcal C}}
\nc{\calE}{{\mathcal E}}
\nc{\calF}{{\mathcal F}}
\nc{\calG}{{\mathcal G}}
\nc{\calH}{{\mathcal H}}
\nc{\calI}{{\mathcal I}}
\nc{\calJ}{{\mathcal J}}
\nc{\calK}{{\mathcal K}}
\nc{\calL}{{\mathcal L}}
\nc{\calM}{{\mathcal M}}
\nc{\calN}{{\mathcal N}}
\nc{\calO}{{\mathcal O}}
\nc{\calP}{{\mathcal P}}
\nc{\calQ}{{\mathcal Q}}
\nc{\calR}{{\mathcal R}}
\nc{\calS}{{\mathcal S}}
\nc{\calT}{{\mathcal T}}
\nc{\calU}{{\mathcal U}}
\nc{\calV}{{\mathcal V}}
\nc{\calW}{{\mathcal W}}
\nc{\calX}{{\mathcal X}}
\nc{\calY}{{\mathcal Y}}
\nc{\calZ}{{\mathcal Z}}

\nc{\Sht}{{\on{Sht}}}
\nc{\Frob}{{\on{Frob}}}
\nc{\Hecke}{{\on{Hecke}}}
\nc{\inv}{{\on{inv}}}
\nc{\Conv}{{\on{Conv}}}
\nc{\triv}{{\on{triv}}}
\nc{\tw}{{\on{tw}}} 
\nc{\Isom}{{\on{Isom}}}

\nc{\scrB}{{\mathscr{B}}}
\nc{\scrA}{{\mathscr{A}}}
\nc{\bbf}{{\mathbf{f}}}
\nc{\bba}{{\mathbf{a}}}
\nc{\rig}{{\mathrm rig}}

\nc{\al}{\alpha}
\nc{\be}{\beta}
\nc{\ga}{\gamma}
\nc{\la}{\lambda}
\nc{\qcqs}{{\on{qcqs}}}
\nc{\Bmu}{{\boldsymbol \mu}}

\nc{\pot}[1]{ [\hspace{-0,5mm}[ {#1} ]\hspace{-0,5mm}] }
\nc{\rpot}[1]{ (\hspace{-0,7mm}( {#1} )\hspace{-0,7mm}) }

\nc{\defined}{\hspace{0.1cm}\stackrel{\text{\tiny \rm def}}{=}\hspace{0.1cm}}
\nc{\co}{\colon}
\nc{\specto}{{\leadsto}}

\newcommand{\category}[1]{\mathrm{#1}}
\newcommand{\Cat}{\category{Cat}} 
\newcommand{\cell}{\mathrm{cell}} 
\newcommand{\Fun}{\category{Fun}} 
\newcommand{\AffSch}{\category{AffSch}} 
\newcommand{\Ani}{\category{Ani}} 
\newcommand{\Sm}{\category{Sm}} 
\newcommand{\Sch}{\category{Sch}} 
\newcommand{\Ab}{\category{Ab}} 
\newcommand{\can}{\mathrm{can}} 
\renewcommand{\Pr}{\category{Pr}}
\newcommand{\PrL}{\Pr^\category{L}} 
\newcommand{\PrSt}{\Pr^{\category{St}}} 
\newcommand{\KGL}{\mathrm{KGL}} 
\newcommand{\MGL}{\mathrm{MGL}} 
\newcommand{\Sp}{\category{Sp}} 

\def\Gm{\mathbf {G}_\mathrm m} 
\newcommand{\GmX}[  1]{\mathbb {G}_{\mathrm {m}, #1}} 
\def\red{\mathrm{red}} 

\font\tencyr=wncyr10
\font\sevencyr=wncyr7
\font\fivecyr=wncyr5
\newcommand{\dual}{\vee} 

\newcommand{\colim}{\operatornamewithlimits{colim}} 
\newcommand{\cofib}{\operatorname{cofib}} 
\newcommand{\coker}{\operatorname{coker}} 
\def\id{{\rm id}} 
\def\opp{{\rm op}} 
\def\To#1#2{\mathop{\count0=#1 \loop\ifnum\count0>0 \smash-\mkern-7mu \advance\count0 -1 \repeat \mathord\rightarrow}\limits^{#2}} 
\def\CH{\mathop{\rm CH}\nolimits} 
\def\THH{\mathop{\rm THH}\nolimits} 
\def\Hom{\mathop{\rm Hom}\nolimits} 
\def\Pic{\mathop{\rm Pic}\nolimits} 
\def\triv{\mathop{\rm triv}} 

\def\Gr{\mathop{\rm Gr}\nolimits} 
\def\Sht{\mathop{\rm Sht}\nolimits} 
\def\Sym{\mathop{\rm Sym}\nolimits} 
\def\CAlg{\mathop{\rm CAlg}\nolimits} 
\def\End{\mathop{\rm End}\nolimits} 
\def\Map{\mathop{\rm Map}\nolimits} 
\def\et{\mathrm{\acute et}} 

\definecolor{hellgrau}{RGB}{200,200,200} 
\definecolor{dunkelgrau}{RGB}{160,160,160} 
\definecolor{hellblau}{RGB}{194, 215, 249} %
\definecolor{dunkelblau}{RGB}{68, 128, 226} %

\def\Z{{\mathbb Z}} 
\def\F{{\mathbb F}} %
\def\Fp{\F_p} %
\def\Fpq{\ol {\F}_p} %
\def\Fq{\F_q} %
\def\N{{\mathbb N}} 
\def\Q{{\mathbb Q}} 
\def\SS{{\mathbb S}} 
\def\Qp{\Q_p} 
\def\A{{\mathbb A}} 
\renewcommand{\P}[1][1]{\mathbb P^{#1}} 
\def\Gm{\mathbb {G}_\mathrm m} 

\def\H{{\rm H}} 
\def\SH{\category{SH}} %
\def\DM{\category{DM}} 
\def\eff{\mathrm{eff}} 
\def\ii{$\infty$}

\def\RG{\R \Gamma} 

\def\Spec{\mathop{\rm Spec}} 
\newcommand{\M}{\mathrm{M}} 
\newcommand{\D}{\category{D}} 

\def\R{{\rm R}} 
\def\sbuildrel#1\over#2{\mathrel{\smash{\mathop{\kern0pt #2}\limits^{#1}}}}

\let\x\times
\let\ol\overline
\renewcommand{\t}{\otimes}

\renewcommand{\r}{\rightarrow}
\newcommand{\lr}{\longrightarrow}

\def\matrix#1{\null\,\vcenter{\normalbaselines
    \ialign{\hfil$##$\hfil&&\quad\hfil$##$\hfil\crcr
      \mathstrut\crcr\noalign{\kern-\baselineskip}
      #1\crcr\mathstrut\crcr\noalign{\kern-\baselineskip}}}\,}

\newdimen\harrowsize
\def\mapright#1{\smash{\mathop{\hbox to\harrowsize{\rightarrowfill}}\limits^{#1}}}

{\catcode`@=11
\gdef\cal{\fam\tw@}
\global\let\over\@@over
\global\let\atop\@@atop
\global\let\above\@@above
\global\let\overwithdelims\@@overwithdelims
\global\let\atopwithdelims\@@atopwithdelims
\global\let\abovewithdelims\@@abovewithdelims
\gdef\eqalign#1{\null\,\vcenter{\openup\jot\m@th
\ialign{\strut\hfil$\displaystyle{##}$&$\displaystyle{{}##}$\hfil
      \crcr#1\crcr}}\,}
\newskip\xcentering \global\xcentering=0pt plus 1000pt minus 1000pt
\gdef\eqalignno#1{\displ@y \tabskip\xcentering
  \halign to\displaywidth{\hfil$\@lign\displaystyle{##}$\tabskip\z@skip
    &$\@lign\displaystyle{{}##}$\hfil\tabskip\xcentering
    &\llap{$\@lign##$}\tabskip\z@skip\crcr
    #1\crcr}}
\global\def\cases#1{\left\{\,\vcenter{\normalbaselines\m@th
    \ialign{$##\hfil$&\quad##\hfil\crcr#1\crcr}}\right.}
\gdef\eqlabel#1{\refstepcounter{equation}\label{eqn--#1}\eqno\hbox{\@eqnnum}}
}

\def \nts#1{}

\def \journal#1
{ 
\noindent\colorbox{dunkelblau}{\parbox{\dimexpr\textwidth-2\fboxsep\relax}{#1}}
}

\begin{document}

\title{Frobenius rigidity in $\A^1$-homotopy theory}

\author{Timo Richarz, Jakob Scholbach}

\thanks{
The first named author T.R.~is funded by the European Research Council (ERC) under Horizon Europe (grant agreement nº 101040935), by the Deutsche Forschungsgemeinschaft (DFG, German Research Foundation) TRR 326 \textit{Geometry and Arithmetic of Uniformized Structures}, project number 444845124 and the LOEWE professorship in Algebra, project number LOEWE/4b//519/05/01.002(0004)/87.
The second named author J.S. acknowledges supported by the European Union – Project 20222B24AY (subject area: PE -
Physical Sciences and Engineering) ``The arithmetic of motives and L-functions'', and logistical support by the Max-Planck-Institute for Mathematics in Bonn.}


\begin{abstract}
We study the homotopy fixed points under the Frobenius endomorphism on the stable $\A^1$-homotopy category of schemes in characteristic $p>0$ and prove a rigidity result for cellular objects in these categories after inverting $p$.
As a consequence we determine the analogous fixed points on the $K$-theory of algebraically closed fields in positive characteristic.
We also prove a rigidity result for the homotopy fixed points of the partial Frobenius pullback on motivic cohomology groups in weights at most $1$.
\end{abstract}

\maketitle

\setcounter{tocdepth}{1}
\tableofcontents

\section{Introduction}
A functor $\kappa \colon \CAlg_k \to \Sp$ from the category of commutative algebras over a field $k$ to the \ii-category of spectra (or the derived category of abelian groups, or groups, or sets) is called \emph{rigid} if for any extension $F \subset E$ of algebraically closed overfields of $k$, the induced map
$$\kappa(F) \r \kappa(E)$$
is an isomorphism.
For example, for some proper $k$-scheme $X$ and some prime number $\ell$, the functor $R \mapsto \H^n_\et(X \x_k \Spec R, \Z/\ell)$ given by taking étale cohomology is rigid for all $n\in \Z$.
Similarly, Suslin's celebrated rigidity result \cite{Suslin:K-theory} states that mod-$\ell$ $K$-theory $R\mapsto K_n(R)/\ell$ is rigid, provided $\ell$ is prime to the characteristic of $k$.
In particular, when $F$ is of characteristic $p>0$ prime to $\ell$, Quillen's computation of $K_n(\Fpq)$ gives explicit results for the mod-$\ell$ $K$-groups of $F$.
Suslin's argument is robust enough to allow for various extensions, including the rigidity result of Röndigs--{\O}stv{\ae}r  \cite{RoendigsOstvaer:Rigidity} asserting the full faithfulness of the pullback functor $\SH(F)/\ell \r \SH(E)/\ell$ between the mod-$\ell$ stable $\A^1$-homotopy categories, for two algebraically closed fields $F \subset E$ of characteristic prime to $\ell$.

Of course, the full stable $\A^1$-homotopy category fails to be rigid, as is visible already for the first $K$-group $K_1(F) = F^\x$.
For $k=\Fp$, the present paper studies the idea of rigidifying various functors by applying (homotopy) fixed points under the Frobenius endomorphism, as opposed to considering classes modulo $\ell$.
As a first indication note that for an algebraically closed field $F$ of characteristic $p>0$ the complex 
$$F^\x \stackrel{x \mapsto x / x^p} \lr F^\x\eqlabel{Fx Frob}$$
is quasi-isomorphic to $\Fp^\x$ by Kummer theory.
 
For an $\Fp$-scheme $S$, let $\Frob_S\colon S\to S$ be the Frobenius endomorphism given by $f\mapsto f^p$ on functions.
If $S$ is understood, we abbreviate $\Frob_S$ simply by $\Frob$.

\defi
For an $\Fp$-scheme $S$, the \emph{Frobenius stable $\A^1$-homotopy category} is the fixed point category
$$\SH(S/\Frob):= \lim \left ( \SH(S)\stackrel[\id]{\Frob^*} \rightrightarrows \SH(S)\right )$$
i.e., the homotopy fixed points of Frobenius acting on the stable $\A^1$-homotopy category.
\xdefi

Objects in this category are pairs of objects $M\in \SH(S)$ together with an isomorphism $M\cong \Frob^*M$.
By construction, the canonical pullback functor $\SH(\Fp)\to \SH(S)$ factors through a functor
$$
\can_S \colon \SH(\F_p)\r \SH(S/\Frob).
$$
The idea of rigidity after taking Frobenius fixed points leads to the following question:

\ques
Is the functor
\eqn\CAlg_{\Fp}\r \Sp,\;R\mapsto \Map_{\SH(R/\Frob)[p^{-1}]}(\can_RM,\can_RN) \label{intro blah}\xeqn
rigid for all $M,N\in \SH(\Fp)[p^{-1}]$?
\xques

The main results of this paper exhibit two situations in which we can answer special cases of this question affirmatively.
To state the first, recall that the subcategory $\SH(\Fp)_\cell \subset \SH(\Fp)$ of \emph{cellular objects} is the full presentable subcategory generated by the motivic spheres $\SS^{m,n}$ for all $m, n \in \Z$.

\theo[\thref{rigidity.cellular}]
The functor in \eqref{intro blah} is rigid for all $M, N \in \SH(\Fp)_\cell[p^{-1}]$.
\xtheo

Applying the theorem to cellular spectra (see also Section \ref{subsec-cellular} for more examples) such as the homotopy invariant $K$-theory spectrum implies the following result, where amusingly $p^{-1}$-localization is not necessary:

\coro[\thref{Frobenius.K}]\thlabel{main:intro}
For any algebraically closed field $F$ of characteristic $p>0$, one has
$$\pi_n(K(F/\Frob)) = \left \{ \begin{tabular}{ll}
$\Z$ & $n = -1, 0$ \\
$\F_{p^i}^\x$ & $n = 2i-1> 0$ \\
$0$ & else, \end{tabular} \right.$$
where $K(F/\Frob)$ denotes the homotopy fixed points of the Frobenius endomorphism on the $K$-theory spectrum $K(F)$.
\xcoro

In the formulation of the next result, we denote by $\SH^\eff(\Fp)$ the stable, full subcategory of $\SH(\Fp)$ generated under colimits by motives of smooth $\Fp$-schemes, but not allowing negative Tate twists.
 
\theo[\thref{Frobenius eff rigid}]
\thlabel{wt 1 intro}
The functor in \eqref{intro blah} is rigid for all $M \in \SH^\eff(\Fp)[p^{-1}]$ and $N=\SS^{n,n}[p^{-1}]$ (or, $N=\Z[p^{-1}](n)$) for all $n \le 1$.
\xtheo

Let us unwind the meaning of this assertion in terms of Bloch's cycle complex.
For a smooth $\Fp$-scheme $X$ of finite type, we define the \emph{Frobenius--Bloch cycle complex}
$$\RG(X \x F / \Frob_F, \Z(n)) := \operatorname{Tot} \left [\RG(X \x F, \Z(n)) \stackrel{\id - (\id_X \x \Frob_F)^*} \lr \RG(X \x F, \Z(n)) \right ]\eqlabel{Tot}$$ 
to be the total complex associated to the two-term double complex, where $\RG(X \x F, \Z(n))$ denotes Bloch's cycle complex.
Equivalently, this is the homotopy fixed point of the action of the partial Frobenius pullback $(\id_X \x \Frob_F)^*$ on $\RG(X \x F, \Z(n))$.
If we take $F = \Fpq$, and consider étale motivic cohomology $\RG_\et(-, \Z(n))$, this recovers Weil-étale cohomology of schemes in characteristic $p>0$ as introduced by Lichtenbaum \cite{Lichtenbaum:Weil} and studied in particular by Geisser \cite{Geisser:Weil}.
We refer to the above concept as \emph{Frobenius motivic cohomology} (as opposed to Weil motivic cohomology) in order to emphasize that fixed points under partial Frobenius are considered even for transcendental extensions $F$ over $\Fp$.
For $M$ being the motive of $X$ and $N = \Z(n)$, the rigidity asserted above amounts to the claim that the complex $\RG(X \x F / \Frob_F, \Z(n))$ is rigid after inverting $p$, i.e., is independent, up to quasi-isomorphism, of the choice of an algebraically closed field $F$ of characteristic $p>0$.

The proof for $n = 1$ is based on the following observation.
Resolution of singularities (by alterations) allows to reduce to the case of $X$ being smooth and proper over $\Fp$.
The maps
\begin{align*}\Gm(X) \x \Gm(Y) \stackrel \boxtimes \r & \Gm(X \x_F Y) \\
\Pic(X) \x \Pic(Y) \stackrel \boxtimes \r & \Pic(X \x_F Y)
\end{align*}
fail to be isomorphisms for algebraically closed fields $F$ in general.
However, the ``error terms'' are under control, cf.~\refeq{Rosenlicht} and \refeq{Pic product}, and the homotopy fixed points of the action by a partial Frobenius on these error terms do vanish.
From this perspective, \thref{wt 1 intro} shares a kinship with a statement known as \emph{Drinfeld's lemma} \cite[Proposition~1.1]{Drinfeld:FSheaves}, which rectifies the failure of étale fundamental groups of $\Fp$-schemes to satisfy a Künneth formula, and is a key point in the Langlands program over fields such as $\Fp(t)$ or $\Qp$ \cite{Lafforgue:Chtoucas,FarguesScholze:Geometrization}.

One may ask whether Frobenius motivic cohomology is rigid for $n \ge 2$ as well.
In that direction, we recall the following result, which is also in the vicinity of Drinfeld's lemma \cite[Lemma~4.7]{HemoRicharzScholbach:Categorical}: for a finite type $\Fp$-scheme $X$, and any  extension of algebraically closed fields $F \subset E$ in characteristic $p>0$, the base change
$$\{\text{constructible subsets of } X \x F \} \r \{\text{constructible subsets of } X \x E \}$$ induces a bijection after restricting to those subsets that are set-theoretically stable under $\id_X \x \Frob_E$, resp.~$\id_X \x \Frob_F$.
In fact, these are precisely the subsets descending to $X$.
This result gives control over the degree-wise kernel of $(\id_X \x \Frob_F)^*$ on Bloch's cycle complex $\RG(X \x F, \Z(n))$. 
The obstacle towards an analogue of \thref{wt 1 intro} for $n \ge 2$ is a similar control of the cokernel.

We conclude this paper with a short appendix on the homotopy fixed points of the partial Frobenius on topological Hochschild homology.
Frobenius THH is again rigid (cf.~\thref{Frobenius THH}), with the Artin--Schreier sequence playing the rôle of the Kummer sequence in the context of Frobenius $K$-theory.

\medskip

\noindent \textbf{Acknowledgements.} 
We thank Tom Bachmann, Markus Land, Zhouhang Mao, Jakob Stix and Georg Tamme for helpful email exchanges and comments on the manuscript.

\section{Rigid functors}
Let $\AffSch_k$ be the category of affine schemes over a field $k$.
We identify its opposite category $\AffSch_k^\opp$ with the category of commutative $k$-algebras $\CAlg_k$ whenever convenient. 
Let $\Ani$ be the \ii-category of anima (also called spaces or \ii-groupoids). 

\defi
\thlabel{rigid-functor}
A functor $\kappa\colon  \AffSch_k^\opp\to \mathcal \Ani$ is called \textit{rigid} if, for any extension $F\subset E$ of algebraically closed fields over $k$, the map $\kappa(F)\to \kappa(E)$ is an equivalence.
\xdefi

Recall that a functor is called {\it finitary} if it preserves filtered colimits.

\lemm[Suslin]
\thlabel{rigidity-criterion}
Let $\kappa\colon \AffSch_k^\opp\to \mathcal \Ani$ be a finitary functor. 
Then, the following are equivalent:
\begin{enumerate}
\item \label{rigidity-criterion:fields}
	The functor $\kappa$ is rigid. 
\item \label{rigidity-criterion:curves}
	For any algebraically closed field $F$ over $k$, any connected, smooth, affine $F$-curve $C$, any $n\in\Z_{\geq 0}$ and any $\alpha \in \pi_n(\kappa(C))$ there exists a non-empty open (automatically affine) subset $U_\alpha\subset C$ such that the map $U_\alpha(E)\to \pi_n(\kappa(E)), c\mapsto c^*\alpha$ is constant for any algebraically closed field extension $E$ over $F$.
\end{enumerate}
\xlemm
\pf
Since equivalences in $\Ani$ are detected on homotopy groups and their formation commutes with filtered colimits, we may and do assume that $\kappa$ takes values in the category of sets.

Let $F$ be an algebraically closed field over $k$.
Then, the map $\kappa(F)\to \kappa(E)$ is injective for any $F$-algebra $E$: 
by finitariness of $\kappa$, and expressing $E$ as the filtered colimit of the finitely generated $F$-subalgebras, we may assume $E$ is a finitely generated $F$-algebra.
Since $F$ is algebraically closed, Hilbert's Nullstellensatz supplies a section of the structural map $F \to E$ implying the injectivity of $\kappa(F)\to \kappa(E)$.
This uses neither \eqref{rigidity-criterion:fields} nor \eqref{rigidity-criterion:curves}.

Now assume \eqref{rigidity-criterion:fields} holds.
Let $C\to \Spec F$ and $\alpha\in \kappa(C)$ be as in \eqref{rigidity-criterion:curves}. 
Let $K$ be an algebraic closure of the function field of $C$.
Then, $K=\colim_{\tilde C\to C} \Gamma(\tilde C,\calO)$ is a filtered colimit where $\tilde C$ ranges over the connected, smooth, affine $F$-curves equipped with a flat (necessarily generically finite) map to $C$. 
Using that $\kappa$ is finitary, we get maps of sets
\[
\kappa(F) \r \kappa(C)\r \kappa(K)=\colim_{\tilde C\to C}\kappa(\tilde C),
\]
whose composition is bijective by \eqref{rigidity-criterion:fields}.
Thus, there exists some $\tilde C\to C$ such that the pullback $\alpha|_{\tilde C}$ lies in the image of $\kappa(F)\to \kappa(\tilde C)$.
In particular, the map $\tilde C(E)\to \kappa(E), \tilde c\mapsto \tilde{c}^*(\alpha|_{\tilde C})$ is constant for any $F$-algebra $E$ where $\tilde C(E)$ denotes the set of $F$-maps $\Spec E\to \tilde C$.
Let $U_\alpha$ be the (necessarily open by flatness) image of $\tilde C\to C$. 
Then, $\tilde{c}^*(\alpha|_{\tilde C})=c^*\alpha$ for $\tilde C(E)\r U_\alpha(E), \tilde c\mapsto c$.
So, \eqref{rigidity-criterion:curves} follows from the surjectivity of $\tilde C(E)\twoheadrightarrow U_\alpha(E)$ for algebraically closed fields $E$.

Conversely, assume that \eqref{rigidity-criterion:curves} holds.
Let $F\subset E$ be an algebraically closed field extension.
It remains to show that the injection $\kappa(F)\hookrightarrow \kappa(E)$ is surjective.
By finitariness of $\kappa$, we reduce to fields $E$ of finite transcendence degree over $F$, then to transcendence degree $1$ by induction.
Again, by finitariness of $\kappa$, any element $\alpha\in\kappa(E)$ arises by pullback from some $\alpha_C\in \kappa(C)$ for a connected, smooth, affine $F$-curve $C$ whose algebraically closed function field identifies with $E$.
Denote by $\eta\in C(E)$ the canonical map.
Using \eqref{rigidity-criterion:curves} and replacing $C$ by $U_\alpha$ if necessary, we may and do assume that the map $C(E)\to \kappa(E), c\mapsto c^*\alpha_C$ is constant.
Applying this to $c=\eta$ gives $\eta^*\alpha_C=\alpha$ by construction.
Hence, choosing any section $\Spec F\r C$ and looking at the composition $\Spec E\to \Spec F\to C$ implies \eqref{rigidity-criterion:fields}.
\xpf

In practice the following corollary is useful:

\coro
\thlabel{geometric-rigidity}
Let $\kappa \colon  \category{AffSch}_k^\opp \r \Ani$ be a finitary functor such that for any algebraically closed field $F$ over $k$, any connected, smooth, affine $F$-curve $C$ and any points $c_0, c_1 \in C(F)$, the maps
  $$\pi_n (\kappa(c_i)) \colon  \pi_n (\kappa(C)) \r \pi_n(\kappa(F)), \;\; i=1,2$$
agree for all $n \in \Z_{\geq 0}$.
Then, $\kappa$ is rigid.
\xcoro
\pf
\thref{rigidity-criterion} \eqref{rigidity-criterion:curves} is satisfied with $U_\alpha=C$ for all $\alpha\in \pi_n (\kappa(C))$, noting that $F$-maps $\Spec E\to C$ are the same as sections of the base change $C\x_F E\to \Spec E$.
\xpf

\rema
\thref{rigid-functor}, \thref{rigidity-criterion} and \thref{geometric-rigidity} admit obvious analogues for finitary functors $\kappa\colon  \AffSch_k^\opp\to \mathcal \Sp$ with values in the \ii-category of spectra $\Sp$, i.e., the stabilization of $\Ani$.
Indeed, in the proof of \thref{rigidity-criterion} one uses the non-degeneracy of the $t$-structure on $\Sp$ and the commutation of $\pi_n$ with filtered colimits to reduce to functors valued in the category of abelian groups $\Ab$.
The rest of the argument is the same.
Likewise, for finitary functors $\kappa\colon  \AffSch_k^\opp\to \mathcal \D(\Z)$ valued in the derived category of abelian groups. 
\xrema

\exam \thlabel{rigid-examples}
The following (non-)examples are of interest throughout:
\begin{enumerate}
	\item 
		Let $M,N\in \SH(k)$ be motivic spectra, see Section \ref{Frobenius-stable:notation} for the definition of $\SH$. 
		Then, for any $n\in \Z$ prime to the characteristic of $k$, the functor $\CAlg_k\to \Sp$ given for any $k$-algebra $R$ by the cofiber 
		\[
		\cofib\left(\Map_{\SH(R)}(M_R,N_R)\overset{n\cdot }{\lr}\Map_{\SH(R)}(M_R,N_R)\right)
		\]
		is rigid \cite{RoendigsOstvaer:Rigidity}. 
		The functor is finitary if $M$ is compact.
	\item\label{rigid-examples:Ktheory}
		$K$-theory defines a finitary functor $K\colon \CAlg_k\to \Sp$ that is not rigid: $K_1(F)=F^\x$ highly depends on $F$.
		In particular, the functor $\Gm\colon \CAlg_k\to \Ab, R\mapsto R^\x$ is not rigid, evidently. 
		This plays well with the fact that for $C = \A^1_F-\{0\}$ and a point $c \in C(F) = F^\x$, the map
		$$\Gm (C) = (F[t^\pm])^\x = F^\x \x t^\Z \stackrel{t \mapsto c} \lr F^\x$$
		depends on the choice of $c$.
		In the subsequent sections, we show that the homotopy fixed points under the (partial) Frobenius endomorphism are rigid. 
\end{enumerate}
\xexam

\section{Fixed point categories}\label{sec:fixedpoints}

\defi
For an endofunctor of an \ii-category $\varphi \colon  \calC \r \calC$ the \emph{fixed point category of $\varphi$ on $\calC$} is the \ii-category
$$\calC^\varphi := \lim \left ( \calC \stackrel[\id]{\varphi} \rightrightarrows \calC \right ) = \calC \x_{\varphi \x \id, \calC \x \calC, \Delta} \calC.$$
\xdefi

\rema
Objects in $\calC^\varphi$ are triples $(c_1 \in \calC, c_2 \in \calC, (c_1, c_1)\cong (\varphi c_2, c_2))$.
Any such object is isomorphic to one of the form $(c, c, (c \cong \varphi c, \id_c))$, i.e., one can think of objects as pairs $(c, c\cong \varphi c)$.
The anima (or space) of maps between two such objects $(c, \lambda\colon c \cong \varphi c)$ and $(c',   \lambda'\colon c' \cong \varphi c')$ is the equalizer in $\Ani$ of the following two maps
$$\xymatrix{
\Map_\calC(c, c') \ar@/^1pc/[rrr]^{\lambda'_*} \ar[r]_\varphi & 
\Map_\calC(\varphi c, \varphi c') \ar[rr]_{\lambda^*} & & \Map_\calC(c, \varphi c').}\eqlabel{equalizer}$$
Thus, maps in $\calC^\varphi$ are maps $f \colon c \r c'$ in $\calC$ together with a commutative diagram:
$$\xymatrix{
c \ar[d]^{f} \ar[r]_\cong^\lambda & \varphi c \ar[d]^{\varphi f} \\
c' \ar[r]_\cong^{\lambda'} & \varphi c'
}\eqlabel{maps}$$
\xrema

\exam
\thlabel{alpha.canonical}
\thlabel{trivial}
For $\varphi = \id_\calC$, one has $\Map_{\calC^\id}(\triv c, \triv c') = \Map_\calC(c,c') \oplus \Map_\calC(c, c')[-1]$, since the two maps $\lambda^* \circ \varphi$ and $\lambda'_*$ in \refeq{equalizer} agree.
In addition, there is a functor
$$\triv \colon \calC \r \calC^\id, c \mapsto (c, c \stackrel \id \r c),$$
which is not fully faithful due to the shifted copy of the mapping spectrum.  
\xexam

\rema
\thlabel{twisting}
The fixed point category $\calC^\varphi$ can also be regarded as the limit of the functor $\Phi_{\calC, \varphi} \colon B \N \r \Cat_\infty$ sending $* \mapsto \calC$ and $\N \ni 1 \mapsto \varphi$. 
In particular, an equivalence of functors $\beta \colon \varphi \r \varphi'$ gives rise to an isomorphism of diagrams $\Phi_{\calC, \varphi} \cong \Phi_{\calC, \varphi'}$, and therefore an equivalence
\begin{equation*}
\calC^\varphi \stackrel[\cong]\beta \r \calC^{\varphi'}, \;\; (c, c \stackrel[\cong] \lambda \r \varphi c) \mapsto (c, c \stackrel [\cong] \lambda \r \varphi c \stackrel[\cong]{\beta (c)} \r \varphi'c).\end{equation*}
So, given an equivalence $\beta\colon \id \stackrel \cong \r \varphi$, 
there is a ``twisting'' functor
$$\tw := \tw_\beta\colon \calC \stackrel \triv \r \calC^\id \stackrel \cong \r \calC^\varphi.$$
\xrema

\rema
\thlabel{Fix symmetric monoidal}
If $\calC$ is a presentably symmetric monoidal (i.e., presentable, symmetric monoidal and the $\otimes$-product commutes with colimits in each variable), stable \ii-category and $\varphi$ a symmetric monoidal endofunctor, then so is $\calC^\varphi$. 
Indeed, the forgetful functors
$$\CAlg(\PrSt) \r \PrSt \r \PrL \r \Cat_\infty.$$
preserve limits, see \cite[Proposition~3.2.2.1]{Lurie:HA} for the first, 
then \cite[Proposition~4.8.2.18]{Lurie:HA}, 
and finally \cite[Proposition~5.5.3.13]{Lurie:Higher}.
In addition, if $\calC$ is compactly generated and $\varphi$ preserves compact objects, then $\calC^\varphi$ is compactly generated \cite[Proof of Lemma 2.5]{HemoRicharzScholbach:Categorical}.
In the situation of \thref{trivial} (or \thref{twisting}, where $\beta$ is supposed to be an equivalence of symmetric monoidal colimit-preserving functors) the functors $\triv$ (resp.~$\tw$) will again be functors in $\CAlg(\PrSt)$.
\xrema

\section{Frobenius stable homotopy category}\label{Frobenius-stable:notation}

For a scheme $S$, we denote by $\SH(S)$ the \textit{stable $\A^1$-homotopy category}, i.e., the presentably symmetric monoidal $\infty$-category given by the $\P$-stabilization of $\A^1$-invariant Nisnevich $\infty$-sheaves of spectra on the category $\Sm_S$ of smooth $S$-schemes, cf.~the discussion around \cite[Equation (C.11)]{Hoyois:Quadratic}, and also \cite[Appendix A]{BachmannHoyois:Norms} for the definition of the Nisnevich topology in full generality.
The construction of $\SH$ gives a functor
$$M\colon \Sm_S\to \SH(S),\eqlabel{taking-the-motive}$$
which associates to a smooth $S$-scheme $X$ its \textit{motive} $M(X)$.
If $S$ is quasi-compact and quasi-separated (qcqs), then $\SH(S)$ is compactly generated \cite[Proposition C.12]{Hoyois:Quadratic}, up to desuspensions, by the motives $M(X)$ of \textit{finitely presented}, smooth $S$-schemes $X$. 
If $S$ is qcqs of finite Krull dimension, then every Nisnevich sheaf is a hypersheaf \cite[Theorem 1.7]{ClausenMathew:Hyper}, so the above definition of $\SH(S)$ agrees with more classical definitions using model categorical language \cite[\S 2.4.1]{Robalo:K-theory}.
The construction of the stable $\A^1$-homotopy category can be upgraded to a functor $\SH\colon \Sch_S^\opp\to \CAlg(\PrSt)$ using $*$-pullback functoriality and further to a six functors formalism \cite{Khan:Motivic, Hoyois:Six}.

We use the following standard notation for the \emph{motivic spheres}: let $\SS^{1,1} \in \SH(S)$ be the object represented by $\GmX S$, and denote by $\SS^{1,0} \in \SH(S)$ the suspension of the monoidal unit.
By definition of $\SH(S)$, both objects are dualizable. 
So, the definition
$$\SS^{n+r,n} := (\SS^{1,1})^{\t n} \t (\SS^{1,0})^{\t r}$$
makes sense for all $n, r \in \Z$.
Note that $\SS^{0,0} = 1$ is the monoidal unit in $\SH(S)$.

\subsection{The Frobenius stable homotopy category}
Fix a prime number $p$.
For an $\Fp$-scheme $S$, we denote by $\Frob_S\colon S\to S$ the Frobenius endomorphism given by $f\mapsto f^p$ on functions. 
If $S$ is understood, we abbreviate $\Frob_S$ simply by $\Frob$.
The pullback $\Frob^*$ induces a symmetric monoidal endofunctor of $\SH(S)$, so the setting of Section \ref{sec:fixedpoints}, in particular \thref{Fix symmetric monoidal}, applies.

\defi
The \emph{Frobenius stable $\A^1$-homotopy category of $S$} is the fixed point category under the pullback along the Frobenius map:
$$\SH(S / \Frob) := \SH(S)^{\Frob^*}$$
\xdefi

\rema\thlabel{variants-fixedpoint-cats}
We use an analogous notation also for other \ii-categories:
\begin{enumerate}
  \item 
  If $P$ is a set of prime numbers, then we denote by $\SH(S)[P^{-1}]$ the full subcategory in $\SH(S)$ of $P^{-1}$-localized objects, i.e., $M\in \SH(S)$ with $M\otimes 1/\ell= 0$ for all $\ell\in P$.
  The inclusion is right adjoint to the localization functor $\SH(S)\to \SH(S)[P^{-1}]$, see e.g.~\cite[\S 3.2]{MathewNaumannNoel:Descent} for a general discussion. 
  Further, $P^{-1}$-localization commutes with taking Frobenius fixed points, and we denote by $\SH(S/\Frob)[P^{-1}]$ the resulting full subcategory of $\SH(S/\Frob)$.
  We apply this to the cases where $P=\{p\}$ and where $P$ is the set of all primes.
  The resulting categories of $p^{-1}$-localized and rational objects are denoted by $\SH(S/\Frob)[p^{-1}]$ and $\SH(S/\Frob)_\Q$ respectively. 
  \item
  Similarly, we consider the category $\DM(S / \Frob)$, where $\DM$ denotes the category of Beilinson motives with rational coefficients \cite{CisinskiDeglise:Triangulated}.
  \end{enumerate} 
\xrema


\rema
\thlabel{pullback}
The formation of $\SH(S / \Frob)$ is functorial since the Frobenius endomorphism is so.
That is, for a map $s \colon  S' \r S$, the $*$-pullback induces a symmetric monoidal functor
\begin{align*}
s^* \colon  \SH( S/\Frob) \r & \SH( S' / \Frob), \\
(M, M \stackrel[\cong]{\lambda} \r \Frob^* M) \mapsto & (s^* M, s^* M \stackrel[\cong]{s^* \lambda} \r s^* \Frob^* M = \Frob^* s^* M).
\end{align*}
\xrema

\cons\thlabel{canonical} \thref{pullback} applies to the structural map $s \colon  S \r \Spec \Fp$ and gives the symmetric monoidal functor
\begin{align*}
\can_S \colon  \SH(\Fp) & \stackrel \triv \lr \SH( \Fp / \id) \stackrel{s^*} \lr \SH(S / \Frob), \\
M \mapsto & (s^* M,s^*M\overset{\id}{\to} s^*M=\Frob^*s^*M).
\end{align*}
using $\Frob=\id$ over $\Fp$ and so $s\circ \Frob=s$.
We also use the same notation for the variants in \thref{variants-fixedpoint-cats}.
\xcons

\subsection{Twisted Frobenius objects}
If $t \colon  T \r S$ is a morphism of $\Fp$-schemes, we  consider throughout the usual diagram involving the relative Frobenius $\Frob_{T/S}$ where the square is cartesian:
$$\xymatrix{
T \ar[r]^{\Frob_{T/S}} \ar[dr]_t & T' \ar[d]^{t'} \ar[r] & T \ar[d]^t \\
& S \ar[r]^{\Frob_S} & S.}\eqlabel{relative Frobenius}$$

\exam
\thlabel{Frob.Gm}
The relative Frobenius $\Frob_{\GmX S/S} = \Frob_{{\mathbb G}_{{\rm m},\Fp}} \x \id_S$ agrees with the $p$-multiplication of the $S$-group scheme $\GmX S$.
\footnote{$$(-)^p \colon  \Gm \stackrel{\Delta} \lr \Gm \x \dots \x \Gm \stackrel{\mu} \lr \Gm$$
is given on the level of Hopf algebras by
$$\Fp[t^{\pm 1}] \stackrel {t \mapsto t \t \dots \t t} \lr \bigotimes_{p \text{ times}} \Fp[t^{\pm1}] \stackrel{\mathrm{mult}} \lr \Fp[t^{\pm1}],$$
which satisfies $t \mapsto t^p$.}
\xexam

\prop
\thlabel{id.Frob}
Let $S$ be an $\Fp$-scheme.
Then, there is a canonical isomorphism of symmetric monoidal endofunctors on $\SH(S)[p^{-1}]$,
$$\beta \colon  \id \stackrel \cong \r \Frob^*$$
given on motives of smooth $S$-schemes $T$ by the relative Frobenius maps:
$$\beta (\M(T)) \colon  \M(T) \stackrel{\Frob_{T/S}} \lr \M(T \x_{S, \Frob} S) = \Frob^* \M(T).$$
\xprop

\pf
In order to construct $\beta$
we use the universal property of $\SH$, see \cite{Robalo:K-theory}. 
Using the notation in \refeq{relative Frobenius} 
 gives a functor
$$\tw \colon  \Sm_{ S} \r \Fun(\Delta^1, \Sm_{ S}), T \mapsto (T \r T')$$
whose evaluations at the two endpoints of $\Delta^1$ are the identity, respectively $\Frob^*$.
This functor has a symmetric monoidal structure with respect to the pointwise monoidal structure on the target category.
By the universal property of $\SH$ \cite[Corollary~2.39]{Robalo:K-theory}, it descends to a symmetric monoidal functor
$$\SH( S)[p^{-1}] \r \Fun(\Delta^1, \SH( S)[p^{-1}])$$
whose evaluations at the two endpoints of $\Delta^1$ are again $\id$ and $\Frob^*$:
On the unstable $\A^1$-homotopy category $\H(S)$ (i.e., on $\A^1$-invariant Nisnevich sheaves), this gives a map $\id \r \Frob^*$ without inverting $p$. 
Its evaluation at the object represented by $\GmX S$ is the $p$-multiplication (\thref{Frob.Gm}), and therefore we obtain $\beta$ after inverting $p$.

It remains to show that $\beta(M)$ is an equivalence for all $M \in \SH( S)$.
It suffices to do this for $M = \M_{ S}(T) = t_! t^! 1_S$ for some smooth $t \colon  T \r  S$ as above.
The map $\Frob_{T/S}$ is a universal homeomorphism \StP{0CCB}. 
So, the functor $\Frob_{T/S}^* \colon  \SH(T')[p^{-1}] \r \SH(T)[p^{-1}]$ is an equivalence \cite{ElmantoKhan:Perfection} with inverse $\Frob_{T/S,*} = \Frob_{T/S,!}$, hence $\Frob_{T/S}^* = \Frob_{T/S}^!$ as well.
We have $t = t' \circ \Frob_{T/S}$ with notation as in \refeq{relative Frobenius}.
So, the counit map $\Frob_{T/S,!}\Frob_{T/S}^!\to \id$, which is an isomorphism, induces the isomorphism  
$$\M_{S}(T) = t_! t^! 1_S = t'_! \Frob_{T/S, !} \Frob_{T/S}^! t'^! 1_S \stackrel \cong \r t'_! t'^! 1_S = \M_{ S}(T') = \Frob^* \M_{ S}(T).$$
It agrees with the map induced by $\Frob_{T/S}\colon T\to T'$ under the functor \refeq{taking-the-motive} on motives.
\xpf

Applying \thref{twisting} to \thref{id.Frob}, we get: 
\coro
For any $\Fp$-scheme $S$, there is the symmetric monoidal ``\emph{twisting}'' functor
$$\tw\colon  \SH(S)[p^{-1}] \stackrel \triv \r \SH(S / \id)[p^{-1}] \stackrel [\cong]{\beta} \r \SH(S / \Frob)[p^{-1}].$$  
For a smooth $\Fp$-scheme $X$, one has
$$\tw(\M(X)) = (\M(X \x S), \M(X \x S) \stackrel{\Frob_X \x \id_S} \lr \M(X \x S)).$$
\xcoro

\rema
\thlabel{two.functors}
The two functors
$$\SH(\Fp)[p^{-1}] \stackrel[\tw]{\triv} \rightrightarrows \SH( \Fp / \id)[p^{-1}]$$
do not agree. Indeed, $\tw(\SS^{r+n,n}) = (\SS^{r+n,n}, p^n \cdot \id)$, by \thref{Frob.Gm}.
\xrema

\rema
\thlabel{pullback.twist}
For a map $s\colon  S' \r S$, the pullback functor from \thref{pullback} together with the twisting functors give a diagram
$$\xymatrix{
\SH(S)[p^{-1}] \ar[d]^{s^*} \ar[rr]^{\tw} & &
\SH(S/\Frob)[p^{-1}]\ar[d]^{s^*} \\
\SH(S')[p^{-1}] \ar[rr]^{\tw} &&  \SH(S'/\Frob)[p^{-1}],
}$$
which commutes since forming relative Frobenii is functorial.
\xrema

\section{Frobenius rigidity}\label{Frobenius-section}

Recall the functor $\can_S\colon \SH(\Fp)[p^{-1}]\to \SH(S/\Frob)[p^{-1}]$ from \thref{canonical}.

\defi
\thlabel{Frobenius rigidity}
An ordered pair of objects $M, N \in \SH(\Fp)[p^{-1}]$ is \emph{Frobenius rigid} if the functor $\AffSch_{\Fp}^\opp\to \Sp$ given by
\[
S\mapsto \Map_{\SH(S / \Frob)[p^{-1}]}(\can_S M, \can_S N)
\]
is rigid.
That is, if for any extension of algebraically closed fields $f\colon \Spec E\to \Spec F$ in characteristic $p>0$ the induced map
$$\Map_{\SH(F / \Frob)[p^{-1}]}(\can_F M, \can_F N) \stackrel {f^*}\lr \Map_{\SH( E / \Frob)[p^{-1}]}(\can_E M, \can_E N)\eqlabel{Frobenius rigidity}$$
is an equivalence. 
\xdefi

\rema
This definition suggests the question to what extent the functor 
$$f^* \colon  \SH(F/\Frob)[p^{-1}] \r \SH(E / \Frob)[p^{-1}]$$
is fully faithful. 
On the whole of $\SH(F/\Frob)[p^{-1}]$, $f^*$ is \emph{not} fully faithful, however.
Indeed, using the twisting functor $\tw$, both categories are equivalent to $\SH(-)[p^{-1}]^{\id}$.
By \thref{trivial} (and given that the functor $f^*$ then identifies with the usual $f^*$, by \thref{pullback.twist}), $f^*$ is not fully faithful, compare also \thref{rigid-examples}~\eqref{rigid-examples:Ktheory}.
\xrema


\rema
If the pair $M, N$ is Frobenius rigid, then the invariance of $\SH(-)[p^{-1}]$ under perfection \cite{ElmantoKhan:Perfection} implies a similar rigidity property for any extension of separably (as opposed to algebraically) closed fields.
\nts{For conditions of closure cf. Poonen Rational points on varieties, Ex. 1.1. https://math.mit.edu/~poonen/papers/Qpoints.pdf}
\xrema

Recall the variants of the Frobenius fixed point categories from \thref{variants-fixedpoint-cats}.

\lemm
\thlabel{Z.Q}
Let $M, N \in \SH(\Fp)[p^{-1}]$.
The following are equivalent:
\begin{enumerate}
  \item \label{item--SH} The pair $M, N$ is Frobenius rigid.
  \item \label{item--SHQ} Their rationalizations $M_\Q, N_\Q$ satisfy the property of \refeq{Frobenius rigidity} in $\SH(- / \Frob)_\Q$.
  \item \label{item--DMQ} The Beilinson motives associated with $M_\Q, N_\Q$ satisfy the property of \refeq{Frobenius rigidity} in $\DM(- / \Frob)$.
\end{enumerate}
\xlemm

\pf
Let $A$ be the fiber of the map in \refeq{Frobenius rigidity}.
By definition, $p$-multiplication is invertible.
The arithmetic fracture square \cite[(3.17)]{MathewNaumannNoel:Descent} 
implies that $A = 0$ if and only if both its rationalization $A_\Q = 0$ and $A / n := \cofib (A \stackrel {n \cdot} \lr A) = 0$ for all $n$ prime to $p$.
Röndigs--{\O}stv{\ae}r's version of Suslin rigidity for $\SH$, i.e., the full faithfulness of $\SH(F)/n \r \SH(E)/n$ \cite{RoendigsOstvaer:Rigidity}, ensures that the latter holds for any $M, N$ as above.
This proves \refit{SH} $\Leftrightarrow$ \refit{SHQ}.

For any field containing a square root of $-1$, in particular for $K=E$ and $K=F$, $\SH(K)_\Q = \DM(K)$ \cite[Corollary~16.2.14]{CisinskiDeglise:Triangulated}.
Again passing to homotopy fixed points under Frobenius pullback shows the equivalence of \refit{SHQ} and \refit{DMQ}.
\xpf

\subsection{Frobenius rigidity for cellular objects}

Recall, e.g.~from \cite[§2.8]{DuggerIsaksen:Motivic} that the subcategory
$$\SH(S)_\cell \subset \SH(S)$$
of \emph{cellular objects} 
is the stable full subcategory generated under colimits by the spheres $\SS^{r+n,n}$, which lie in the essential image of $\SH(\Fp)\to \SH(S)$, for all $r, n \in \Z$.
These objects are dualizable and, if $S$ is qcqs, also compact.

\theo
\thlabel{rigidity.cellular}
Any pair of $p^{-1}$-localized cellular objects $M, N \in \SH(\Fp)_\cell[p^{-1}]$ is Frobenius rigid.
\xtheo

\pf
As $\SH(\Fp)_\cell$ is compactly generated by dualizable objects, we may assume $M = 1_{\Fp}$.
Again using compactness, we may then also assume that $N = \SS^{r+n, n}$ is a compact generator of $\SH(\Fp)_\cell$.
Let $F\subset E$ be an extension of algebraically closed fields in characteristic $p>0$.
We have to prove that the map between the $p^{-1}$-localized mapping spectra 
$$\Map_{\SH(F / \Frob)}(1, \can_F \SS^{r+n,n})[p^{-1}] \r \Map_{\SH(E / \Frob)}(1, \can_E \SS^{r+n,n})[p^{-1}]\eqlabel{Map.E.F}$$ 
is an isomorphism for all $r, n \in \Z$.
Since $\SS^{1,0}$ is the object associated with the constant presheaf with values the circle, we may assume $r=0$.

Letting $S$ denote either $\Spec F$ or $\Spec E$, we will show that these mapping spectra are insensitive to the choice of $\Spec F$ or $\Spec E$.
By definition, $\can_S1_{\Fp} = (1_S, \can_{1_S} \colon  1_S \cong \Frob_S^* 1_S)$, see \thref{canonical}.
Abbreviating $\SS := s^*\SS^{n,n}$, the same description holds for $\can_S\SS^{n,n} = (\SS, \can_\SS \colon  \SS \cong \Frob_S^* \SS)$.
We have the following canonical identifications, where all mapping spectra appearing at the right are in $\SH(S)[p^{-1}]$:
\begin{align*}
\Map_{\SH(S / \Frob)[p^{-1}]}(1, \can_S \SS^{n,n}) & = \lim \left ( \Map(1, \SS) \stackrel[\can_\SS \circ -]{\Frob_S^* (-) \circ \can_{1_S}} \rightrightarrows \Map(1, \Frob_S^* \SS) \right ) \\
& = \lim \left ( \Map(1, \SS) \stackrel[\id]{\can_{\SS}^{-1} \circ \Frob_S^* (-) \circ \can_{1_S}} \rightrightarrows \Map(1, \SS) \right ) \\
& = \lim \left ( \Map(1, \SS) \stackrel[\id]{\can_{\SS}^{-1} \circ \beta(\SS) \circ -} \rightrightarrows \Map(1, \SS) \right ) \\
& = \lim \left ( \Map(1, \SS) \stackrel[\id]{p^n \cdot - } \rightrightarrows \Map(1, \SS) \right ).
\end{align*}
The first equality follows from \refeq{equalizer}, the second by postcomposing with $\can_\SS^{-1}$, the equality $\Frob_S^* (-) \circ \can_{1_S}=\beta(\SS)$ from $\can_{1_S} = \beta(1_S)$ and the functoriality of $\beta$ (\thref{id.Frob}) and the last from \thref{Frob.Gm}, according to which the composite $\can_\SS^{-1} \beta(\SS)$ equals $p^n \cdot \id$.
Thus, it remains to show that the fiber of the multiplication with $1-p^n$ on $\Map_{\SH(S)[p^{-1}]}(1, \SS)=\Map_{\SH(S)}(1, \SS)[p^{-1}]$ is insensitive to replacing $S = \Spec F$ by $S=\Spec E$. 

By Suslin rigidity (\thref{Z.Q}), we may consider the category of Beilinson motives $\DM(-)$ instead of $\SH(-)[p^{-1}]$.
Then, each homotopy group of the associated mapping spectra is a $\Q$-vector space, so multiplication by $1-p^n\neq 0$ is an isomorphism for all $n \ne 0$. 
Therefore $\Map_{\DM(S/\Frob)}(1, \can_S \Q(n)[n]) = 0$ in this case.
For $n=0$, already the mapping spectra $\Map_{\DM(S)}(\Q,\Q) = \Q$ 
are independent of the chosen $S$. 
Passing to homotopy fixed points under the trivial Frobenius actions preserves that independence.
\xpf


\subsection{Frobenius stable homotopy groups}

Recall that the \emph{stable $\A^1$-homotopy groups} of a field $F$ are defined as 
$$\pi_{r,n}(F) := \Hom_{\SH(F)}(\SS^{r+n,n}, 1),$$
where $1$ denotes the monoidal unit.

Morel showed that these groups vanish for $r < 0$  \cite[Theorem~4.9]{Morel:A1}.
For $r=0$ they are isomorphic to \emph{Milnor--Witt $K$-groups} $K_{-n}^{MW}(F)$, which for algebraically closed fields reduce to the Milnor $K$-groups $K_{-n}^M(F)$.
For odd primes $p$ and any $p$-power $q$, the $p^{-1}$-localized groups $\pi_{1,n}(\Fq)[p^{-1}]$ have been computed in \cite[\S 8.10]{OrmsbyOstvaer:Stable}.
We also have $K^M_n (\Fq) = 0$ for $n \ge 2$ \cite[Example~1.5]{Milnor:Algebraic}.
These computations, and the continuity of $\SH$, which allows to pass to $\Fpq = \colim \Fq$, give the following results for $\pi_{r,n}(\Fpq)[p^{-1}]$ and odd primes $p$:
\begin{center}
\begin{tabular}{c|ccccccc}
$n$ & $\le-2$ & $-1$ & 0 & 1 & 2 & $\ge 2$ \\
\hline
$r = 1$ & 0 & 0 
& $(\Z/2)^{\oplus 2}$ 
& $\Z/2$ 
& $\Z/24[p^{-1}]$
& 0
\\
$r= 0$ & 0
& $\Fpq^\x$ 
& $\Z[p^{-1}]$ 
& 0 & 0 & 0
\end{tabular}
\end{center}

\defi
Let $F$ be a field of characteristic $p>0$. 
The \emph{$p^{-1}$-localized Frobenius stable $\A^1$-homotopy groups} are the groups
$$\pi_{r,n}(F / \Frob)[p^{-1}] := \Hom_{\SH(F / \Frob)[p^{-1}]}(\can_F \SS^{r+n, n}, 1).$$
\xdefi

These groups appear in a long exact sequence:
$$\dots \r \pi_{r,n}(F / \Frob)[p^{-1}] \r \pi_{r,n}(F)[p^{-1}] \stackrel{\id - \Frob} \lr \pi_{r,n}(F)[p^{-1}] \lr \pi_{r-1,n}(F / \Frob)[p^{-1}] \r \dots$$
The Frobenius rigidity of cellular spectra (\thref{rigidity.cellular}) implies the following computation:

\coro
\thlabel{Frobenius.Milnor.K}
The groups $\pi_{r,n}(F / \Frob)[p^{-1}]$ are independent of the choice of an algebraically closed field $F$ of characteristic $p>0$.
For small values of $r$, and odd primes $p$, the groups are given by
\begin{center}
\begin{tabular}{c|ccccccc}
$n$ & $\le-2$ & $-1$ & 0 & 1 & 2 & $\ge 2$ \\
\hline
$r = 0$ & $0$ & $\Fp^\x$ & $(\Z/2)^{\oplus 2} \oplus \Z[p^{-1}]$ 
& $\Z/2$
& $\Z/24[p^{-1}]$ & $0$ 
\\
$r= -1$ & $0$
& $0$ 
& $\Z[p^{-1}]$ 
& $0$ & $0$ & $0$
\end{tabular}
\end{center}
\xcoro

\subsection{Frobenius $K$-theory}

\defi
The \emph{Frobenius $K$-theory spectrum} of $S$, with respect to an $\Fp$-scheme $X$, is defined as the equalizer in the \ii-category of spectra
$$K(X \x S/\Frob_S) := \lim \left ( K(X \x S) \stackrel[\id]{(\id_X \x \Frob_S)^*} \rightrightarrows K(X \x S) \right ),$$
i.e., the homotopy fixed points of the pullback along the partial Frobenius $\id_X \x \Frob_S$. The homotopy groups of this spectrum, denoted by $K_n(X \x S/\Frob_S)$, appear in a long exact sequence
$$\dots \lr K_n(X \x S / \Frob_S) \lr K_n(X \x S) \stackrel{\id - \Frob_S^*} \lr K_n(X \x S) \lr K_{n-1}(X \x S / \Frob_S) \lr \dots .\eqlabel{K.Frobenius}$$
\xdefi

In order to relate these groups to the Frobenius fixed points on $\SH$, we place the spectrum $\KGL \in \SH(\Fp)$ representing homotopy $K$-theory inside $\SH(S/\Frob)$ as follows:

\defi
\thlabel{KGL.Frob}
For a scheme $S$ of characteristic $p>0$, let $\KGL_S/\Frob := \can_S \KGL \in \SH(S / \Frob)$.
\xdefi

By \thref{canonical}, $\KGL_S / \Frob $ consists of the spectrum $\KGL_S$ together with the map
$$\can_{\KGL_S}\colon \KGL_S = s^* \KGL \stackrel{\id} \lr s^* \KGL = \Frob_S^* s^* \KGL,$$
where $s\colon  S\to \Spec \Fp$ denotes the structure map.
Being the image of a commutative algebra object in $\SH(\Fp)$, the object $\KGL_S / \Frob$ again has the structure of a commutative algebra object in $\SH(S/\Frob)$.
This object represents Frobenius $K$-theory as follows:

\lemm
\thlabel{Frobenius.K.SH}
Let $S$ be regular Noetherian, and let $X$ be smooth of finite type over $\Fp$.
Then, there is an isomorphism of spectra
$$K(X \x S / \Frob_S) = \Map_{\SH(S/\Frob)}(\can_S \M(X), \KGL_S/\Frob).$$
\xlemm

\pf
We have $s^* \M(X) = \M(X \x S)$. By the assumptions, $X \x S$ is regular,
 so that there is an identification of mapping spectra $\Map_{\SH(S)}(s^* \M(X), \KGL_S) = \Map_{\SH(X \x S)}(1, \KGL_{X \x S}) = K(X \x S)$.
By construction \cite[\S 13.1]{CisinskiDeglise:Triangulated}, for an endomorphism $\varphi$ of $X \x S$, such as $\varphi = \id_X \x \Frob_S$, the map
$$\Map_{\SH(X \x S)}(1, \KGL_{X \x S}) \stackrel {\varphi^*} \lr \Map(\varphi^* 1, \varphi^* \KGL_{X \x S}) \stackrel[\cong]{\can_\KGL^{-1}\circ(-)\circ\can_1} \r \Map(1, \KGL_{X \x S})$$
identifies with the pullback $\varphi^*$ on the $K$-theory spectrum.
(Here at the right $\can_\bullet$ denotes again the canonical isomorphisms coming from functoriality of *-pullback, see \thref{canonical}).
The following computation is analogous to the proof of \thref{rigidity.cellular}, where $\Map_{-} := \Map_{\SH(-)}$:
\begin{align*}
& \Map_{\SH(S / \Frob)}(\can_S \M(X), \can_S \KGL) \\
& = \lim \left ( \Map_S(\M(X \x S), \KGL_S) \stackrel[\can_{\KGL} \circ -]{\Frob_S^* (-) \circ \can_{\M(X \x S)}} \rightrightarrows \Map_S(\M(X \x S), \Frob_S^* \KGL)  \right ) \\
& = \lim \left ( \Map_S(\M(X \x S), \KGL) \stackrel[\id]{\can_{\KGL}^{-1} \circ \Frob_S^* (-) \circ \can_{\M(X \x S)}} \rightrightarrows \Map_S(\M(X \x S), \KGL) \right ) \\
& = \lim \left ( \Map_{X \x S}(1, \KGL_{X \x S}) \stackrel[\id]{\can_{\KGL}^{-1} \circ (\id_X \x \Frob_S)^* (-) \circ \can_{\M(X \x S)}} \rightrightarrows \Map_{X \x S}(1, \KGL_{X \x S}) \right ) \\
& = \lim \left ( K(X \x S) \stackrel[\id]{(\id_X \x \Frob_S)^*} \rightrightarrows K(X \x S) \right ) =: K(X \x S / \Frob_S).
\end{align*}
\xpf

The following result asserts that the Frobenius acts so richly on the $K$-theory of (large enough) fields that hardly anything is fixed under Frobenius pullback.

\coro
\thlabel{Frobenius.K}
Frobenius $K$-theory is rigid. 
That is, for an extension $F \subset E$ of algebraically closed fields in characteristic $p>0$, the pullback map
$$K(F / \Frob) \r K(E / \Frob)$$
is an equivalence of spectra. 
The individual Frobenius $K$-groups are given by
$$K_n(F / \Frob) = \left \{ 
\begin{tabular}{ll}
$\Z$ & $n = -1,0$ \\
$\F_{p^i}^\x$ & $n = 2i-1 > 0$ \\
$0$ & else.
\end{tabular} \right .  \eqlabel{Frobenius.K.groups}$$
\xcoro

\pf
The statement is clear for $n \le 0$, cf.~the discussion around \refeq{Fx Frob}.
By \cite[Theorem~5.4]{Hiller:lambda-rings}, the groups $K_n(F)$ are uniquely $p$-divisible for $n > 0$.
Thus $K_n(F / \Frob) = K_n(F/\Frob)[p^{-1}]$ for $n > 0$.
The first statement now follows from the cellularity of $\KGL$ \cite[Theorem~6.2]{DuggerIsaksen:Motivic}, \thref{Frobenius.K.SH} for $X=\Spec \Fp$ and \thref{rigidity.cellular}.

To see the concrete values in \refeq{Frobenius.K.groups}, we may assume $F=\Fpq$ and use Quillen's computation of $K(\Fpq)$ and its Frobenius action, as reported e.g., in \cite[§VI.1, p.~465]{Weibel:Kbook}: as an abelian group $K_{2i-1}(\Fpq)$ is isomorphic to $\Fpq^\x$, with $\Frob^*$ acting by raising to the $p^i$-th power.
Then, our statement follows from the Kummer sequence.
\nts{$$1 \lr \F_{p^i}^\x \lr \Fpq^\x \stackrel{x \mapsto x^{p^i-1}} \lr \Fpq^\x \lr 1.$$}
\xpf

\rema
As communicated to us by Georg Tamme,
\thref{Frobenius.K} can be proven directly by using that for $n > 0$, the map $\Frob^*$ on $K_n(F)$ agrees with the $p$-th Adams operation, which acts by multiplication with $p^k$ on the $k$-th Adams eigenspace inside $K_n(F)_\Q$.
Such an argument seems not applicable to cellular objects other than $\KGL$.
\xrema

\coro
\thlabel{Frobenius.K.any.field}
The rationalized Frobenius $K$-groups of \emph{any} field $F$ of characteristic $p>0$ are given by
$$K_n(F / \Frob)_\Q = \left \{ 
\begin{tabular}{ll}
$\Q$ & $n = -1, 0$ \\
0 & else.
\end{tabular} \right .  \eqlabel{Frobenius.K.groups.Q}$$
In particular, the \emph{Beilinson--Soulé vanishing} holds for Frobenius $K$-theory of fields:
$$\H^p(F / \Frob, \Q(q)) = K_{2q-p}(F / \Frob)_\Q^{(q)} = 0$$
for $q > 0$ and $p \le 0$.
\xcoro

\pf
The cases $n=-1, 0$ are clear. 
Suppose now $n \ne -1,0$. 
To show the claimed vanishing, we may assume $F$ is perfect, since $p^{-1}$-localized $K$-theory is insensitive to perfection.
Let $\overline F$ be an algebraic closure of $F$.
Combining \thref{Frobenius.K} with finitaryness of Frobenius $K$-theory we have
\begin{equation}
0 = K_n(\overline F / \Frob)_\Q = \colim_{F \subset L \subset \overline F} K_n(L / \Frob)_\Q,\label{colimit.extension}
\end{equation}
where the colimit runs over the finite, separable extensions $F \subset L \subset \overline F$.

It suffices to show that the transition maps $f^* \colon  K(F/\Frob)_\Q \r K(L/\Frob)_\Q$ in \eqref{colimit.extension} are injective.
By \thref{basic functoriality} below, the usual $f_* f^*$ on K-theory extends to Frobenius K-theory.
The maps $f_* f^* \colon  K_n(F) \r K_n(F)$ are equal to $[L:F] \cdot \id$ \cite[§7, Proposition 4.8]{Quillen:HAK}, and are therefore isomorphisms after passing to rationalizations.
So, $f^* \colon  K(F/\Frob)_\Q \r K(L/\Frob)_\Q$ is injective.
\xpf

\lemm
\thlabel{basic functoriality}
If $f \colon  S' \r S$ is a finite étale map, there is a natural pushforward map $f_* \colon  K(X \x S'/\Frob_{S'}) \r K(X \x S/\Frob_S)$, compatible with the usual pushforward on K-theory. 
The same holds for pullback along arbitrary maps $f$.
\xlemm

\pf
The map $f^*$ always exists since $f^* \Frob_S^* = \Frob_{S'}^* f^* \colon  K(X \x S) \r K(X \x S')$.
For étale maps $f$, the natural map $S' \r S' \x_{f, S, \Frob_S} S$ is an isomorphism \StP{0EBS}, so that the base-change formula
$$\Frob_S^* f_* = f_* \Frob_{S'}^* \colon  K(X \x S') \r K(X \x S)$$
implies the existence of the pushforward on Frobenius K-theory.
\xpf

Recall a conjecture of Beilinson: for all fields $F / \Fp$ the canonical map
$$K^M_*(F)_\Q \to K_*(F)_\Q$$
is an isomorphism. 
This conjecture is implied by the Bass conjecture or, alternatively, also by the Tate conjecture, see \cite[Introduction]{GeisserLevine:K-theory} for references.
The next result confirms this conjecture for the Frobenius variants of these two theories.
For an abelian group $A$, we write $A_{(p)} := A \t \Z_{(p)}$ for the localization at the prime ideal $(p)$.

\coro
For any field $F$ of characteristic $p>0$, and any $n \in \Z$, the map 
$$[K_n^M(F)_{(p)} \stackrel{\id - \Frob} \lr K_n^M(F)_{(p)}] \r [K_n(F)_{(p)} \stackrel{\id - \Frob} \lr K_n(F)_{(p)}]$$
is a quasi-isomorphism.
\xcoro

\pf
By \cite[Proof of Theorem~8.1]{GeisserLevine:K-theory}, the natural map $K^M_n(F) \r K_n(F)$ is injective and its cokernel $C$ is a $\Z[p^{-1}]$-module.
We therefore have $C_{(p)} = C_\Q$.
Thus, we may replace the localization at $(p)$ in the statement by the rationalization.

If $\sqrt{-1} \in F$, then $K^M_n(F)_\Q = 0$ for $n > 0$ \cite[Theorem~1.4]{Milnor:Algebraic}. Thus, the claim for algebraically closed fields follows immediately from \thref{Frobenius.Milnor.K} and \thref{Frobenius.K}.

For arbitrary $F$, the argument from the proof of \thref{Frobenius.K.any.field} carries over:
Milnor $K$-theory is continuous, and for any finite field extension $E \subset F$, the composite $K^M_*(E) \r K^M_*(F) \r K^M_*(E)$ is $[F : E] \cdot \id$.
(This is one of the joint properties of Milnor $K$-theory and $K$-theory, cf.~also \cite[Axiom~R2d]{Rost:Chow}.)
\xpf

\rema
The localization at the prime ideal $(p)$ is necessary in the statement above:
for $n > 0$, the group $K_n^M(\F_2)$ vanishes, but $K_{2i-1}(\F_2) = \F_{2^i}^\x \ne 0$ \cite[Corollary~IV.1.13]{Weibel:Kbook}.
\xrema

\subsection{Further cellular spectra}\label{subsec-cellular}

In addition to 1 and $\KGL$, further cellular spectra include the cobordism spectrum $\MGL$ \cite[Theorem~6.4]{DuggerIsaksen:Motivic} as well as, for $p \ne 2$, hermitian $K$-theory and Witt theory \cite[Theorem~1.1]{RoendigsSpitzweckOstvaer:Cellularity}.
Thus, \thref{Frobenius.K} admits analogues for Frobenius cobordism groups, Frobenius hermitian K- and Witt groups.

\subsection{Motivic cohomology of small weight}

An interesting case of \thref{Frobenius rigidity} not covered by the results thus far is the case $M = \M(X)[p^{-1}]$ and $N = \SS^{n,n}[p^{-1}]$ (or $N = \Z[p^{-1}](n)$ or $\Q(n)$, which makes no difference in view of \thref{Z.Q}).
For an $\Fp$-scheme $X$ of finite type and any field $F$ of characteristic $p>0$, we write
$\RG(X \x F, \Z(n))$ for Bloch's complex of codimension $n$-cycles on $X\x F$.
Its $m$-th homology identifies with the higher Chow group $\CH^n(X\x F, m)$. 
As before, we define the Frobenius variant of this by taking homotopy fixed points under the partial Frobenius pullback:
$$\RG(X \x F / \Frob_F, \Z(n)) := \lim \left (\RG(X \x F, \Z(n)) \stackrel[\id]{(\id_X \x \Frob_F)^*} \rightrightarrows \RG(X \x F, \Z(n)) \right )$$
A concrete representative for this complex is the total complex of a two-step double complex, as in \refeq{Tot}.
The cohomology groups of this complex, denoted by $\H^*(X \x \Spec F / \Frob_F, \Z(n))$, again sit in a long exact sequence similar to \refeq{K.Frobenius}.

\theo
\thlabel{rigidity weight 1}
Let $X$ be smooth, proper $\Fp$-scheme.
Then, Frobenius rigidity holds for the pair $\M(X)[p^{-1}]$ and $\SS^{n,n}[p^{-1}]$ for all $n \le 1$.
In particular, the $p^{-1}$-localized Frobenius motivic cohomology groups
$$\H^*(X \x \Spec F / \Frob_F, \Z(n))[p^{-1}]\eqlabel{Frobenius.motivic}$$
are independent of the choice of an algebraically closed field $F$ of characteristic $p>0$, for all $n \le 1$.
\xtheo

After some preparation, the proof will be given at the end of this subsection.

\exam
\timo{added}
For $F=\Fpq$, the étale version of \refeq{Frobenius.motivic} is studied in \cite{Lichtenbaum:Weil, Geisser:Weil} and in \cite{HemoRicharzScholbach:Categorical, HemoRicharzScholbach:Constructible} for constructible $\ell$-adic sheaves.
\xexam

\subsubsection{Rigidity of Frobenius units}
Starting with the non-rigid presheaf $\Gm$ (\thref{rigid-examples}), we do get a rigid functor once we apply homotopy fixed points under the partial Frobenius:

\lemm
\thlabel{units.rigid}
For a geometrically connected and geometrically reduced scheme $X$, the following functor is rigid:
$$\Gm (X \x - / \Frob_{-}) \colon  \CAlg_{\Fp} \r \D(\Z),\; R \mapsto \left [ \Gm(X \x R) \stackrel{\id - \Frob_R^*} \lr \Gm(X \x R) \right ].\eqlabel{unit_rigid}$$
\xlemm
\pf
Since the functor is finitary it suffices to check the criterion in \thref{geometric-rigidity}.
Let $C = \Spec R$ be a connected, smooth, affine curve over an algebraically closed field $F$.
We first consider the case when $X$ is geometrically integral. We can then apply the unit theorem due to Sweedler \cite{Sweedler:Units} and Rosenlicht \cite{Conrad:Units} to $X_F$ and $C$ over $F$ and obtain the following short exact sequences:
$$\xymatrix{
1 \ar[r] & F^\x \ar[d]^{\id-\Frob_F^*} \ar[r] & \Gm(X_F) \x \Gm(C) \ar[d]^{\id-(\Frob_F^*,\Frob_C^*)} \ar[r] &\Gm(X \x C) \ar[d]^{\id-\Frob_C^*} \ar[r] & 1 \\
1 \ar[r] & F^\x \ar[r] &\Gm(X_F) \x \Gm(C) \ar[r] &\Gm(X \x C) \ar[r] & 1,
}\eqlabel{Rosenlicht}$$
which are compatible with the displayed vertical maps.

Below, we write $\mathrm{(co)ker}$ for the (co)kernel of the vertical maps.
Using the snake lemma along with the Kummer sequence, we obtain a short exact sequence and an isomorphism:
$$1 \r \Fp^\x \r \ker|_{\Gm(X_F) \oplus \Gm(C)} \r \ker|_{\Gm(X \x C)} \r 1$$
$$\coker|_{\Gm(X_F) \oplus \Gm(C)} \stackrel \cong \r \coker|_{\Gm(X \x C)}.$$
For $f\in \Gm(X \x C)$ and $c\in C(F)$, the pullback $c^*f\in \Gm(X_F)$ is clearly independent of $c$ if $f$ factors over the projection $X \x C \r X_F$.
Thus, the rigidity of $R \mapsto \text{(co)ker}|_{\Gm(X \x R)}$ follows from the one of $R\mapsto \text{(co)ker}|_{\Gm(R)}$, i.e., we may and do assume $X = \Spec \Fp$.
Then, the pullback map $c^* \colon  \text{(co)ker}|_{\Gm(C)} \r \text{(co)ker}|_{\Gm(F)}$ is independent of $c$:

For the kernel, the left hand group is $\Fp^\x$ since $C$ is integral, and the value of constant functions is clearly independent of $c$.
For the cokernel, it is independent since the target group $\coker|_{\Gm(F)}$ is trivial because $F$ is algebraically closed.
\xpf

\coro\thlabel{rational.units.rigid}
For a geometrically reduced $\Fp$-scheme $X$ with finitely many geometric connected components, the rationalization $\Gm (X \x - / \Frob_{-})_\Q \colon  \CAlg_{\Fp} \r \D(\Q)$ of \refeq{unit_rigid} is rigid.
\xcoro
\pf
First off, we have $\pi_0(X_{\Fpq})=\pi_0(X_F)$ for any algebraically closed field $F$ of characteristic $p$, see \StP{0363}.
Clearly, the rows in diagram \refeq{Rosenlicht} remain exact when replacing $F^\x$ at the left by $G:=\Z[\pi_0(X_F)] \t_\Z F^\x$. 
On this group, $\Frob_F^*$ acts as usual on $F^\x$ and by permutation on the set $\pi_0(X_F)$.
We claim that, after rationalization, the map $\id-\Frob_F^*$ is invertible on $G_\Q$.
In particular, its cokernel vanishes and the proof of \thref{units.rigid} carries over.
To show the claim, let $\Z$ act on $G_\Q$ through $1\mapsto \Frob_F^*$.
We have to show that $\H^i(\Z,G_\Q)=0$ for $i=0$ (injectivity) and $i=1$ (surjectivity).
Observe that $(\Frob_F^*)^n$ acts as the identity on the finite set $\pi_0(X_F)$ for some suitable $n\in\Z_{\geq 1}$.
Thus, the Kummer sequence shows $\H^i(n\Z,G_\Q)=0$ for $i=0,1$.
This obviously implies $\H^0(\Z,G_\Q)=0$. 
The vanishing of $\H^1$ now follows from the inflation-restriction exact sequence 
$$0 \r \H^1(\Z/n\Z, \H^0(n\Z, G_\Q)) \r \H^1(\Z, G_\Q) \r \H^1(n\Z, G_\Q).$$

\xpf

\subsubsection{Verschiebung}

To show rigidity of a Frobenius version of the Picard group, we use some generalities about the \emph{Verschiebung} of abelian varieties, see, e.g., \cite[§5.2]{EdixhovenGeerMoonen:Abelian}.
\nts{
Let $A / S$ be a commutative flat group scheme over a scheme $S$ in characteristic $p$.
The diagonal map $\Delta \colon  A \r A^p$ and the sum map $m \colon  A^p \r A$ are both $\Sigma_p$-equivariant (since $A$ is commutative), and there is a commutative diagram
$$\xymatrix{
A \ar@{=}[r] \ar[d]^\Delta & A \ar[d] \ar[drr]^{\Frob_A} \ar[dr]|{\Frob_{A/S}} \\
\H_0(\Sigma_p, A^p) = \Sym^p A \ar[d]^m & \hat \H_0(\Sigma_p, A^p) \ar[d]^m \ar[l] \ar[r]^{\varphi_{A/S}} & A^{(p)} \ar@{.>}[dl]^{V_{A/S}} \ar[d] \ar[r] & A \ar[d] \\
A \ar@{=}[r] & A & S \ar[r]^{\Frob_S} & S.
}$$
The notation $\H_0$ refers to group homology, i.e., coinvariants (in the category of schemes) with respect to the $\Sigma_p$-action.
The notation $\hat \H_0$ refers to Tate group homology: locally, if $A = \Spec R \r S = \Spec B$ is affine, then $\H_0(\Sigma_p, A^p) = \Spec (R \t_B \dots \t_B R)^{\Sigma_p} = \Spec (\H^0(\Sigma_p, R^{\t_B p})$ and $\hat \H_0(\Sigma_p, A^p) = \Spec \hat \H^0(\Sigma_p, R^{\t_B p})$, i.e., the spectrum of invariants modulo norms.
The relative Frobenius map $\Frob_{A/S}$ factors as indicated.
The flatness of $A/S$ implies that $\varphi_{A/S}$ is an isomorphism \cite[Lemma~5.17]{EdixhovenGeerMoonen:Abelian}.
This allows to define the Verschiebung as
$$V_{A/S} \colon  A^{(p)} \stackrel{m \circ \varphi_{A/S}^{-1}} \lr A.$$
The entire diagram is functorial with respect to maps of commutative flat $S$-group schemes $A \r A'$.
In particular, the Verschiebung commutes with any morphism of such group schemes.}
Recall that for an abelian variety $A$ over $\Fp$, the Verschiebung is an isogeny
$$V_A \colon  A \lr A$$
of degree $p^{\dim A}$.
It commutes with any morphism of abelian varieties $A \r A'$.

Fix a (geometrically) normal, proper $\Fp$-scheme $X$.
We consider the Verschiebung of the Picard variety $A := \Pic^0_{X / \Fp,\red}$, which is an abelian variety over $\Fp$. 
Indeed, $\Pic_{X/\Fp}^0$ is a geometrically irreducible, proper $\F_p$-group scheme \cite[Lemma~9.5.1, Theorem~9.5.4, Remark~9.5.6]{Fantechi:FGAexplained}. 
Its reduction $\Pic_{X/\Fp,\red}^0$ is geometrically reduced and still an $\Fp$-group scheme (since $\Fp$ is perfect both properties are clear, but also hold over general fields by \cite[Discussion above Theorem~5.1.1]{ColliotTheleneSkorobogatov:Brauer}), hence an abelian variety \StP{03RO}.  



\lemm
\thlabel{Frob.Verschiebung}
The Verschiebung of $A=\Pic^0_{X / \Fp,\red}$ and the map induced by pulling back line bundles along the Frobenius $\Frob_X$ agree:
$$V_A = \Frob_X^*.$$
\xlemm

\pf
It suffices to see 
$$\Frob_X^* \circ \Frob_A = V_A \circ \Frob_{A},$$ 
because $\Frob_A$ is an epimorphism (since $A$ is reduced).
By construction of the Verschiebung, the composition $V_A \circ \Frob_A$ is multiplication by $p$.

The simple, but crucial observation (e.g., \cite[Lemme~1.4]{Stroh:Parametrisation}) is that the map $\Frob_A$ sends a $T$-point $a\colon  T \r A$ to $\Frob_A \circ a = a \circ \Frob_T$.
Interpreting $a$ as a line bundle $\calL$ on $X \x_{\Fp} T$, this means that $\Frob_A(\calL) = (\id_X \x \Frob_T)^* \calL$.
Composing this with $\Frob_X^*$, we see that it gets sent to $(\Frob_X \x \Frob_T)^* \calL = \Frob_{X \x T}^* \calL$.
Generally, pulling back line bundles along the total Frobenius on a scheme, such as $X \x T$, sends $\calL$ to $\calL^{\t p}$, as can be seen by regarding the transition functions, which are raised to their $p$-th power. 
Hence, $\Frob_X^* \circ \Frob_A$ is also the $p$-multiplication.
\xpf

\prop
\thlabel{Verschiebung.isogeny}
For any abelian variety $A$ over $\Fp$, and any $\lambda \in \Q$, the element
$$\id + \lambda V_A$$
is an isogeny, i.e., an invertible element in $\End(A)_{\Q}$.
\nts{It is questionable whether for $r \in \Z$ and $s \in Mat_{n \x n}(\Z)$ we still have that 
$$r \id - s V_{A^n} \colon  A^n \r A^n$$
is an isogeny. E.g. for $2 \x 2$-matrices, we may get $\sqrt p$ as an eigenvalue?!?}
\xprop

\pf
Using that the Verschiebung is compatible with any morphism of abelian varieties, we may replace $A$ by any isogeneous abelian variety $A'$ to check this claim (since then $\End(A)_\Q = \End(A')_\Q$).
Therefore we may assume $A = \prod A_i$ is a product of simple abelian varieties $A_i / \Fp$.
The morphism $\id + \lambda V_A$ respects this product decomposition, so we may assume $A$ is simple.
Then $\End(A)_\Q$ is a skew field, so it suffices to show that $\id + \lambda V_A$ is a non-zero element in $\End(A)_\Q$.
The case $\lambda = 0$ being trivial, we now consider $\lambda = \frac r s \in \Q$ with $r, s \in \Z \setminus \{0\}$.
If $s \id_A = r V_A$, 
\nts{$\End(A)$ is a free abelian group, \cite[Corollary 1, §16, p. 178]{Mumford:Abelian}.}
then taking degrees ($\deg V_A = p^{\dim A}$, \cite[Proposition~5.20]{EdixhovenGeerMoonen:Abelian}), we get 
$$s^{2 \dim A} = r^{2 \dim A} p^{\dim A},$$
which is a contradiction.
\xpf

\rema
\thlabel{odd.power.remark}
If $q$ is an {\it odd} $p$-power, then the analogue of \thref{Verschiebung.isogeny} holds for abelian varieties $A/\bbF_q$ equipped with their Verschiebung $V_{A/\bbF_q}$.
Indeed, noting that $\deg V_{A/\bbF_q}=q^{\dim A}$ the same arguments lead to the equation $s^{2 \dim A} = r^{2 \dim A} q^{\dim A}$ which contradicts the assumption that $\log_p(q)$ is odd.
 \xrema

\subsubsection{Rigidity for Frobenius--Picard groups}

\prop
\thlabel{Pic.rigid}
Let $X$ be a smooth, proper $\Fp$-scheme.
Then, the following functor is rigid:
$$\Pic (X \x - / \Frob_{-})_\Q \colon  \AffSch_{\Fp}^\opp \r \D(\Q), S \mapsto \left [ \Pic(X \x S)_\Q \stackrel{\id - (\id_X \x \Frob_S)^*} \lr \Pic(X \x S)_\Q \right ].$$
\xprop

\pf
We have $p = \Frob_{X \x S}^* = \Frob_X^* \circ \Frob_S^*$ on $\Pic(X\x S)$.
In particular, $\Frob_X^*$ is invertible on the rationalization $\Pic(X \x S)_\Q$.
The above complex is therefore quasi-isomorphic to $[\Pic(X \x S)_\Q \stackrel{p - \Frob_X^*} \lr \Pic(X \x S)_\Q]$.

Let $F$ be an algebraically closed field of characteristic $p$, $C$ a smooth, affine, connected $F$-curve, and let $c_0, c_1\in C(F)$ be points. 
Write $X_F = X \x F$.
In order to show rigidity, let $\ol C$ be the smooth compactification of $C$, and let $D = \ol C \setminus C$ be the boundary points.
There is an exact sequence
$$0 \r \Gm(X_F \x_F \ol C) \r \Gm(X_F \x_F C) \r \Z^{\pi_0(X\x D)} \r \Pic(X_F \x_F \ol C) \r \Pic(X_F \x_F C) \r 0$$
using that $X_F$ is smooth and proper.
The sequence is functorial under $\Frob_X^*$. 
On the rationalization of the middle term, $\Q^{\pi_0(X\x D)}$, the map $p - \Frob_X^*$ is easily seen to be invertible: $\Frob_X^*$ acts through permutation on the finite set $\pi_0(X\x D)$, so its eigenvalues on $\Q^{\pi_0(X\x D)}$ are roots of unity.

Therefore, we may assume $C$ is projective in the sequel.

We compute the Picard group of $X \x C = X_F \x_F C$ using the short exact sequence \cite[(5.31)]{ColliotTheleneSkorobogatov:Brauer}
$$0 \r \Pic(X_F) \oplus \Pic(C) \r  \Pic(X \x C) \r \Hom_{\category{AbVar}_F}(B^\dual, A) \r 0,\eqlabel{Pic product}$$
where $B = \Pic^0_{C / F}$ is the Picard variety of the smooth, projective, connected curve $C$, $B^\dual$ its dual abelian variety, and $A = \Pic^0_{X_F / F, \red}=\Pic^0_{X/\Fp,\red} \x F$ the Picard variety of $X_F$.
This sequence is compatible with $(\Frob_X \x \id_F)^* \oplus \id$, resp.~$(\Frob_X \x \id_C)^*$, resp. the map $A \r A$ induced by pullback along $\Frob_X$ (and $\id_F$).
By \thref{Frob.Verschiebung}, the map induced by pullback along $\Frob_X$ on the reduced Picard scheme $\Pic^0_{X/\Fp, \red}$ is the Verschiebung; here we use the assumptions on $X$.
In $\End(A)_\Q$, the element $[p]_A - \Frob_X^* = [p]_A - V_A$ is invertible by \thref{Verschiebung.isogeny}, if $\dim A > 0$.
This implies that postcomposing with $[p]_A - \Frob_X^*$ is an isomorphism on $\Hom(B^\dual, A)_\Q$ (if $\dim A = 0$, this Hom-group is trivial).
Thus, the ``error term'' $\Hom(B^\dual, A)_\Q$ vanishes after passing to homotopy fixed points under $p - \Frob_X^*$.
The restriction of $c_i^*$ on the subgroup $\Pic(X_F) \oplus \Pic(C)$ is clearly independent of the point $c_i \in C(F)$, since it is the identity on $\Pic(X_F)$ and $0 \colon  \Pic(C) \r \Pic(F) = 0$.

\xpf

\pf[Proof of \thref{rigidity weight 1}]
The case $n \le 0$ is trivial since $\RG(X \x \Spec F, \Z(n)) = 0$ for $n < 0$ and is quasi-isomorphic to $\Z[0]$ for $n = 0$.  

We now turn to $n = 1$, using that $\Z(1) = \Gm[-1]$. The only non-zero groups $\H^r(X \x \Spec F, \Gm)$ are for $r = 0$ and $r = 1$, so it suffices to show that the two-term complex
$$[\H^r(X \x \Spec F, \Gm) \stackrel{\id - \Frob_F^*} \lr \H^r(X \x \Spec F, \Gm)]$$
is insensitive (up to quasi-isomorphism) to the choice of an algebraically closed field, at least after $p^{-1}$-localization.
By Suslin rigidity (\thref{Z.Q}), it suffices to consider the rationalization of these two-term complexes.

The formation of this complex is finitary in $F$. 
Our claim then follows for $r = 0$ by \thref{rational.units.rigid} and for $r =1 $ by \thref{Pic.rigid}.
\xpf

\rema
It would be interesting to apply the above ideas towards \emph{Gabber rigidity} for Frobenius motivic cohomology, along the lines of \cite[§4]{GilletThomason:K-theory}. More precisely, one can ask whether for a Henselian local ring $A$ of a smooth variety over an algebraically closed field in characteristic $p$, with residue field $k$, the map
$$\H^n(X \x A / \Frob_A, \Z(1)[p^{-1}]) \r \H^n(X \x k / \Frob_k, \Z(1)[p^{-1}])$$
is an isomorphism.
\footnote{
Some comments on this task.
\lemm
Let $\kappa \colon  \CAlg_{\Fp} \r \Ab$ be a continuous functor such that for a henselian ring $A$ and a smooth affine curve $C \r \Spec A$ the map
$$C(A) \x \kappa(C) \r \kappa(A), (x, \beta) \mapsto x^* \beta$$
does not depend on the choice of the point $x \colon  \Spec A \r C$.
Then 
$$\kappa(R) = \kappa(R\{t\})$$
(henselization of $R[t]$ at the origin).

If, moreover, $\kappa$ is invariant under universal homeomorphisms, then it satisfies Gabber rigidity:
$$\kappa(R) = \kappa(k)$$
($R$ any Henselian ring with residue field $k$).

In this case the Nisnevich sheaf associated to $\kappa$, on the category of smooth $\Fp$-schemes, is the constant sheaf associated to $\kappa(\Fp)$.
\xlemm

\pf
This is an abstraction of the proof of \cite[Theorem~4.1, Corollaries~4.2, 4.4]{GilletThomason:K-theory}.
\xpf

This lemma shall be applied to $\Gm(X \x_{\Fp} - / \Frob)$ and similarly Frobenius-Pic.
The key thing is to make sure that 
\refeq{Rosenlicht} and \refeq{Pic product} remain true, where now instead of an algebraically closed field one has $S = \Spec A$ (pro-smooth Henselian local ring), $X$ will be replaced by the pullback $X \x_{\Fp} S$, and $Y$ will be the smooth curve $C / S$.
By the setup of Gillet-Thomason's rigidity lemma, $C \r S$ has a section (and the point is to show independence of the section).
Since $S$ is (check!) over an algebraically closed field, $X \x_{\Fp} S \r S$ also has a section.
The exactness of the sequence \refeq{Rosenlicht} at the two left spots is then formal.
To show exactness of the right, one can use the formulation in equation (1) in \cite{Conrad:Units}.
Note that $S$ is a local domain by the setup. 
Let $\eta = Q(S)$ be its generic point. 
One has to make sure that $\Gm((X \x_{\Fp} S) \x_S C)$ embedds into the group of units on $((X \x_{\Fp} S) \x_S C) \x_S \eta$.
This will follow if the scheme $(X \x_{\Fp} S) \x_S C$ is integral. Ideas for this: $C \r S$ should be fiberwise geometrically connected, and $X / \Fp$ is by assumption geometrically connected, and hence this asserted connectedness.
Once this embedding is OK, one can check the formula (1) there by passing to the generic fiber, where it holds by \cite{Conrad:Units}.

\refeq{Pic product}: one needs to ensure the existence of a smooth projective compactification of the curve $C$, relative over $S = \Spec A$. This seems to be handled in \url{https://arxiv.org/pdf/1809.04158.pdf}, Lemma 8.2 there.
\cite[Proposition~5.7.1]{ColliotTheleneSkorobogatov:Brauer} seems to carry over for $X \x_S Y$, as long as $X \r S$, $Y \r S$ has a section, and if the two maps are fiberwise (over $S$) geometrically integral. These sections exist per the above.
Over a general base scheme the properties of the Picard scheme are reported in Kleiman: 3.7 + 4.8 there: if $X / S$ is projective flat and fiberweise geometrically integral then $\Pic_{X/S}$ exists. If $X /S$ has a section then the conditions in 2.5 are met. 
If $X/S$ is a relative curve (projective, flat with geometrically integral fibers) then $\H^2(\mathcal O_{X_s}) = 0$ and then $\Pic^0_{X/S}$ is an abelian scheme. (Ex. 5.23, Rk 5.26 in Kleiman), i.e., smooth. It is not fully clear how to handle $\Pic_{X \x S / S}$, but it should be just $\Pic_{X/\Fp} \x S$.
One should investigate whether $(\Pic_{X_S / S})_\red$ is isomorphic to $(\Pic_{X/\Fp})_\red \x S$, using here that $S$ is pro-smooth.
Then investigate whether the sequence in (5.31) in \cite{ColliotTheleneSkorobogatov:Brauer} remains true.}
\xrema

To round off the discussion concerning Frobenius motivic cohomology of small weight, we consider the stable, full subcategory $\SH^\eff(\Fp)$ in $\SH(\Fp)$ generated under colimits by motives of smooth $\Fp$-schemes $X$.

\coro
\thlabel{Frobenius eff rigid} 
Suppose $M \in \SH^{\eff}(\Fp)[p^{-1}]$ and $N=\SS^{n,n}[p^{-1}]$ (or $N = \Z[p^{-1}](n)$) with $n \le 1$.
Then, the pair $M, N$ is Frobenius rigid.
\xcoro

\pf
By resolution of singularities (via alterations), it is known that $\SH^\eff(\Fp)[p^{-1}]$ is  the stable, full subcategory generated under colimits by $\M(X)(e)[e]$, with $X / \Fp$ being smooth and proper, and $e \ge 0$ \cite[Theorem~2.4.3]{BondarkoDeglise:Dimensional}.
Thus, the corollary follows from \thref{rigidity weight 1}.
\xpf

\appendix
\section{Frobenius topological Hochschild homology}

In this aside, we consider homotopy fixed points under Frobenius pullbacks for topological Hochschild homology (THH).
Since $\THH$ is not representable in $\SH$, the following result is not strictly an example of Frobenius rigidity as in \thref{Frobenius rigidity}, but may still be illustrational.

We fix an $\Fp$-scheme $X$.
Recall, e.g., from \cite{NikolausScholze:Topological} the topological Hochschild homology functor
$$\THH (X \x -)\colon  \Sch_{\Fp}^\opp \r \Sp.$$
We let Frobenius THH be again the homotopy fixed points of partial Frobenius:
$$\THH (X \x S / \Frob_S) := \lim \left (\THH(X \x S) \stackrel[\id]{(\id_X \x \Frob_S)^*} \rightrightarrows \THH(X \x S) \right ).$$

\prop
\thlabel{Frobenius THH}
Let $X$ be a smooth, affine $\Fp$-scheme.
Then, Frobenius THH with respect to $X$ is rigid. 
More precisely, for any algebraically closed field $F$ of characteristic $p$, the following natural map is an equivalence:
$$\THH(X) \stackrel \cong \r \THH(X \x \Spec F / \Frob_F).$$
\xprop
\pf
By the Hochschild--Kostant--Rosenberg theorem for $\THH$ due to Hesselholt \cite[Theorem~B]{Hesselholt:p-typical},\footnote{or see \cite[Proposition~5.6]{Morrow:Topological} for a recent exposition} there is an isomorphism\nts{This uses that $F$ is perfect.}
\begin{align*}
\THH_n(X \x \Spec F) & = \bigoplus_{i \ge 0} \Omega^{n-2i}_{X \x F/F} = \bigoplus_{i \ge 0} \Omega^{n-2i}_{X / \Fp} \t_{\Fp} F.
\end{align*}
The Artin--Schreier sequence $0 \lr \Fp \lr F \stackrel{x \mapsto x^p-x} \lr F \r 0$ shows that the homotopy fixed points of $(\id_X \x \Frob_F)^*$ acting on this agree with $\THH_n(X)$.\nts{The map is surjective since $F$ is separably closed and the polynomial $t^p - t - a$ is separable; so it has a zero.}
\xpf

\rema
As Markus Land and Zhouhang Mao have independently pointed out, \thref{Frobenius THH} holds for all $\Fp$-schemes $X$ by using the symmetric monoidal structure of $\THH$ to reduce to the case $X=\Spec\Fp$ and then invoking B\"okstedt periodicity.
\xrema


\begin{thebibliography}{MNN17}

\bibitem[BD17]{BondarkoDeglise:Dimensional}
Mikhail Bondarko and Fr\'{e}d\'{e}ric D\'{e}glise.
\newblock Dimensional homotopy t-structures in motivic homotopy theory.
\newblock {\em Adv. Math.}, 311:91--189, 2017.
\newblock \href {https://doi.org/10.1016/j.aim.2017.02.003}
  {\path{doi:10.1016/j.aim.2017.02.003}}.

\bibitem[BH21]{BachmannHoyois:Norms}
Tom Bachmann and Marc Hoyois.
\newblock Norms in motivic homotopy theory.
\newblock {\em Ast\'{e}risque}, (425):ix+207, 2021.
\newblock \href {https://doi.org/10.24033/ast} {\path{doi:10.24033/ast}}.

\bibitem[CD19]{CisinskiDeglise:Triangulated}
Denis-Charles Cisinski and Fr\'{e}d\'{e}ric D\'{e}glise.
\newblock {\em Triangulated categories of mixed motives}.
\newblock Springer Monographs in Mathematics. Springer, 2019.
\newblock \href {https://doi.org/10.1007/978-3-030-33242-6}
  {\path{doi:10.1007/978-3-030-33242-6}}.

\bibitem[CM21]{ClausenMathew:Hyper}
Dustin Clausen and Akhil Mathew.
\newblock Hyperdescent and \'{e}tale {$K$}-theory.
\newblock {\em Invent. Math.}, 225(3):981--1076, 2021.
\newblock \href {https://doi.org/10.1007/s00222-021-01043-3}
  {\path{doi:10.1007/s00222-021-01043-3}}.

\bibitem[Con]{Conrad:Units}
Brian Conrad.
\newblock Units on product varieties.
\newblock URL: \url{http://math.stanford.edu/~conrad/papers/unitthm.pdf}.

\bibitem[CTS21]{ColliotTheleneSkorobogatov:Brauer}
Jean-Louis Colliot-Th{\'e}l{\`e}ne and Alexei~N. Skorobogatov.
\newblock {\em The {Brauer}-{Grothendieck} group}, volume~71 of {\em Ergeb.
  Math. Grenzgeb., 3. Folge}.
\newblock Cham: Springer, 2021.
\newblock \href {https://doi.org/10.1007/978-3-030-74248-5}
  {\path{doi:10.1007/978-3-030-74248-5}}.

\bibitem[DI05]{DuggerIsaksen:Motivic}
Daniel Dugger and Daniel~C. Isaksen.
\newblock Motivic cell structures.
\newblock {\em Algebr. Geom. Topol.}, 5:615--652, 2005.
\newblock URL: \url{https://doi.org/10.2140/agt.2005.5.615}.

\bibitem[Dri87]{Drinfeld:FSheaves}
V.~G. Drinfeld.
\newblock Moduli varieties of {$F$}-sheaves.
\newblock {\em Funktsional. Anal. i Prilozhen.}, 21(2):23--41, 1987.

\bibitem[EK20]{ElmantoKhan:Perfection}
Elden Elmanto and Adeel~A. Khan.
\newblock Perfection in motivic homotopy theory.
\newblock 2020.

\bibitem[EvdGM]{EdixhovenGeerMoonen:Abelian}
Bas Edixhoven, Gerard van~der Geer, and Ben Moonen.
\newblock Abelian varieties.
\newblock URL: \url{http://van-der-geer.nl/~gerard/AV.pdf}.

\bibitem[FGI{\etalchar{+}}05]{Fantechi:FGAexplained}
Barbara Fantechi, Lothar G\"{o}ttsche, Luc Illusie, Steven~L. Kleiman, Nitin
  Nitsure, and Angelo Vistoli.
\newblock {\em Fundamental algebraic geometry}, volume 123 of {\em Mathematical
  Surveys and Monographs}.
\newblock American Mathematical Society, Providence, RI, 2005.
\newblock Grothendieck's FGA explained.
\newblock \href {https://doi.org/10.1090/surv/123}
  {\path{doi:10.1090/surv/123}}.

\bibitem[FS21]{FarguesScholze:Geometrization}
Laurent Fargues and Peter Scholze.
\newblock Geometrization of the local langlands correspondence, 2021.
\newblock \href {http://arxiv.org/abs/2102.13459} {\path{arXiv:2102.13459}}.

\bibitem[Gei04]{Geisser:Weil}
Thomas Geisser.
\newblock Weil-\'etale cohomology over finite fields.
\newblock {\em Math. Ann.}, 330(4):665--692, 2004.
\newblock URL: \url{http://dx.doi.org/10.1007/s00208-004-0564-8}, \href
  {https://doi.org/10.1007/s00208-004-0564-8}
  {\path{doi:10.1007/s00208-004-0564-8}}.

\bibitem[GL00]{GeisserLevine:K-theory}
Thomas Geisser and Marc Levine.
\newblock The {{\(K\)}}-theory of fields in characteristic {{\(p\)}}.
\newblock {\em Invent. Math.}, 139(3):459--493, 2000.
\newblock \href {https://doi.org/10.1007/s002220050014}
  {\path{doi:10.1007/s002220050014}}.

\bibitem[GT84]{GilletThomason:K-theory}
Henri~A. Gillet and Robert~W. Thomason.
\newblock The {K}-theory of strict {Hensel} local rings and a theorem of
  {Suslin}.
\newblock {\em J. Pure Appl. Algebra}, 34:241--254, 1984.
\newblock \href {https://doi.org/10.1016/0022-4049(84)90037-9}
  {\path{doi:10.1016/0022-4049(84)90037-9}}.

\bibitem[Hes96]{Hesselholt:p-typical}
Lars Hesselholt.
\newblock On the {{\(p\)}}-typical curves in quillen's {{\(K\)}}-theory.
\newblock {\em Acta Math.}, 177(1):1--53, 1996.
\newblock \href {https://doi.org/10.1007/BF02392597}
  {\path{doi:10.1007/BF02392597}}.

\bibitem[Hil81]{Hiller:lambda-rings}
Howard~L. Hiller.
\newblock lambda-rings and algebraic {K}-theory.
\newblock {\em J. Pure Appl. Algebra}, 20:241--266, 1981.
\newblock \href {https://doi.org/10.1016/0022-4049(81)90062-1}
  {\path{doi:10.1016/0022-4049(81)90062-1}}.

\bibitem[Hoy14]{Hoyois:Quadratic}
Marc Hoyois.
\newblock A quadratic refinement of the {G}rothendieck-{L}efschetz-{V}erdier
  trace formula.
\newblock {\em Algebr. Geom. Topol.}, 14(6):3603--3658, 2014.
\newblock \href {https://doi.org/10.2140/agt.2014.14.3603}
  {\path{doi:10.2140/agt.2014.14.3603}}.

\bibitem[Hoy17]{Hoyois:Six}
Marc Hoyois.
\newblock The six operations in equivariant motivic homotopy theory.
\newblock {\em Adv. Math.}, 305:197--279, 2017.
\newblock URL: \url{https://doi.org/10.1016/j.aim.2016.09.031}.

\bibitem[HRS23]{HemoRicharzScholbach:Constructible}
Tamir Hemo, Timo Richarz, and Jakob Scholbach.
\newblock Constructible sheaves on schemes.
\newblock {\em Adv. Math.}, 429:46, 2023.
\newblock Id/No 109179.
\newblock \href {https://doi.org/10.1016/j.aim.2023.109179}
  {\path{doi:10.1016/j.aim.2023.109179}}.

\bibitem[HRS24]{HemoRicharzScholbach:Categorical}
Tamir Hemo, Timo Richarz, and Jakob Scholbach.
\newblock A categorical {K}\"{u}nneth formula for constructible {W}eil sheaves.
\newblock {\em Algebra Number Theory}, 18(3):499--536, 2024.
\newblock \href {https://doi.org/10.2140/ant.2024.18.499}
  {\path{doi:10.2140/ant.2024.18.499}}.

\bibitem[Kha16]{Khan:Motivic}
Adeel Khan.
\newblock Motivic homotopy theory in derived algebraic geometry, 2016.
\newblock {PhD thesis, Universit\"at Duisburg-Essen}.

\bibitem[Laf18]{Lafforgue:Chtoucas}
Vincent Lafforgue.
\newblock Chtoucas pour les groupes r\'{e}ductifs et param\'{e}trisation de
  {L}anglands globale.
\newblock {\em J. Amer. Math. Soc.}, 31(3):719--891, 2018.
\newblock URL: \url{https://doi.org/10.1090/jams/897}.

\bibitem[Lic05]{Lichtenbaum:Weil}
S.~Lichtenbaum.
\newblock The {W}eil-\'etale topology on schemes over finite fields.
\newblock {\em Compos. Math.}, 141(3):689--702, 2005.
\newblock \href {https://doi.org/10.1112/S0010437X04001150}
  {\path{doi:10.1112/S0010437X04001150}}.

\bibitem[Lur09]{Lurie:Higher}
Jacob Lurie.
\newblock {\em Higher topos theory}, volume 170 of {\em Annals of Mathematics
  Studies}.
\newblock Princeton University Press, Princeton, NJ, 2009.

\bibitem[Lur17]{Lurie:HA}
Jacob Lurie.
\newblock {Higher Algebra}.
\newblock 2017.
\newblock URL: \url{http://www.math.harvard.edu/~lurie/}.

\bibitem[Mil70]{Milnor:Algebraic}
John~W. Milnor.
\newblock Algebraic {{\(K\)}}-theory and quadratic forms. {With} an appendix by
  {J}. {Tate}.
\newblock {\em Invent. Math.}, 9:318--344, 1970.
\newblock \href {https://doi.org/10.1007/BF01425486}
  {\path{doi:10.1007/BF01425486}}.

\bibitem[MNN17]{MathewNaumannNoel:Descent}
Akhil Mathew, Niko Naumann, and Justin Noel.
\newblock Nilpotence and descent in equivariant stable homotopy theory.
\newblock {\em Adv. Math.}, 305:994--1084, 2017.
\newblock \href {https://doi.org/10.1016/j.aim.2016.09.027}
  {\path{doi:10.1016/j.aim.2016.09.027}}.

\bibitem[Mor06]{Morel:A1}
Fabien Morel.
\newblock {{\(\mathbb{A}^1\)}}-algebraic topology.
\newblock In {\em Proceedings of the international congress of mathematicians
  (ICM), Madrid, Spain, August 22--30, 2006. Volume II: Invited lectures},
  pages 1035--1059. Z{\"u}rich: European Mathematical Society (EMS), 2006.

\bibitem[NS18]{NikolausScholze:Topological}
Thomas Nikolaus and Peter Scholze.
\newblock On topological cyclic homology.
\newblock {\em Acta Math.}, 221(2):203--409, 2018.
\newblock \href {https://doi.org/10.4310/ACTA.2018.v221.n2.a1}
  {\path{doi:10.4310/ACTA.2018.v221.n2.a1}}.

\bibitem[O{\O}14]{OrmsbyOstvaer:Stable}
Kyle~M. Ormsby and Paul~Arne {\O}stv{\ae}r.
\newblock Stable motivic {{\(\pi_1\)}} of low-dimensional fields.
\newblock {\em Adv. Math.}, 265:97--131, 2014.
\newblock \href {https://doi.org/10.1016/j.aim.2014.07.024}
  {\path{doi:10.1016/j.aim.2014.07.024}}.

\bibitem[Qui73]{Quillen:HAK}
Daniel Quillen.
\newblock Higher algebraic {$K$}-theory. {I}.
\newblock In {\em Algebraic {$K$}-theory, {I}: {H}igher {$K$}-theories ({P}roc.
  {C}onf., {B}attelle {M}emorial {I}nst., {S}eattle, {W}ash., 1972)}, pages
  85--147. Lecture Notes in Math., Vol. 341. Springer, Berlin, 1973.

\bibitem[R{\O}08]{RoendigsOstvaer:Rigidity}
Oliver R\"{o}ndigs and Paul~Arne {\O}stv{\ae}r.
\newblock Rigidity in motivic homotopy theory.
\newblock {\em Math. Ann.}, 341(3):651--675, 2008.
\newblock \href {https://doi.org/10.1007/s00208-008-0208-5}
  {\path{doi:10.1007/s00208-008-0208-5}}.

\bibitem[{Rob}15]{Robalo:K-theory}
Marco {Robalo}.
\newblock {$K$-theory and the bridge from motives to noncommutative motives.}
\newblock {\em {Adv. Math.}}, 269:399--550, 2015.
\newblock \href {https://doi.org/10.1016/j.aim.2014.10.011}
  {\path{doi:10.1016/j.aim.2014.10.011}}.

\bibitem[Ros96]{Rost:Chow}
Markus Rost.
\newblock Chow groups with coefficients.
\newblock {\em Doc. Math.}, 1:319--393, 1996.

\bibitem[RS{\O}19]{RoendigsSpitzweckOstvaer:Cellularity}
Oliver R{\"o}ndigs, Markus Spitzweck, and Paul~Arne {\O}stv{\ae}r.
\newblock Cellularity of {Hermitian} {{\(K\)}}-theory and {Witt}-theory.
\newblock In {\em \(K\)-theory. Proceedings of the international colloquium,
  Mumbai, 2016}, pages 35--40. New Delhi: Hindustan Book Agency; Mumbai: Tata
  Institute of Fundamental Research, 2019.

\bibitem[{Sta}17]{StacksProject}
The {Stacks Project Authors}.
\newblock {Stacks Project}.
\newblock \url{http://stacks.math.columbia.edu}, 2017.

\bibitem[Str17]{Stroh:Parametrisation}
Beno\^{\i}t Stroh.
\newblock La param\'{e}trisation de {L}anglands globale sur les corps de
  fonctions.
\newblock {\em Ast\'{e}risque}, (390):Exp. No. 1110, 169--197, 2017.
\newblock S\'{e}minaire Bourbaki. Vol. 2015/2016. Expos\'{e}s 1104--1119.

\bibitem[Sus83]{Suslin:K-theory}
A.~Suslin.
\newblock On the {K}-theory of algebraically closed fields.
\newblock {\em Invent. Math.}, 73:241--245, 1983.
\newblock \href {https://doi.org/10.1007/BF01394024}
  {\path{doi:10.1007/BF01394024}}.

\bibitem[Swe70]{Sweedler:Units}
M.~E. Sweedler.
\newblock A units theorem applied to {Hopf} algebras and {Amitsur} cohomology.
\newblock {\em Am. J. Math.}, 92:259--271, 1970.
\newblock \href {https://doi.org/10.2307/2373506} {\path{doi:10.2307/2373506}}.

\bibitem[Wei13]{Weibel:Kbook}
Charles~A. Weibel.
\newblock {\em The {$K$}-book}, volume 145 of {\em Graduate Studies in
  Mathematics}.
\newblock American Mathematical Society, Providence, RI, 2013.
\newblock An introduction to algebraic $K$-theory.

\end{thebibliography}

\newcommand{\etalchar}[1]{$^{#1}$}

\end{document}